\documentclass[a4paper,12pt]{article}

\makeatletter
\@addtoreset{footnote}{page}
\makeatother

\usepackage{style-enumitem}

\usepackage{amsthm}
\usepackage{amsmath,amssymb,latexsym,amsfonts,mathrsfs}

\newcommand{\n}{\nonumber}

\renewcommand{\hat}{\widehat}

\renewcommand{\bar}{\overline}
\newcommand{\un}[1]{\underline{#1}}

\newcommand{\absol}[1]{\left| #1 \right|} 
\newcommand{\rbra}[1]{\!\left( #1 \right)} 
\newcommand{\cbra}[1]{\!\left\{ #1 \right\}} 
\newcommand{\sbra}[1]{\!\left[ #1 \right]} 

\newcommand{\bE}{\ensuremath{\mathbb{E}}}
\newcommand{\bF}{\ensuremath{\mathbb{F}}}

\newcommand{\bN}{\ensuremath{\mathbb{N}}}

\newcommand{\bP}{\ensuremath{\mathbb{P}}}

\newcommand{\bR}{\ensuremath{\mathbb{R}}}

\newcommand{\cA}{\ensuremath{\mathcal{A}}}
\newcommand{\cB}{\ensuremath{\mathcal{B}}}

\newcommand{\cF}{\ensuremath{\mathcal{F}}}

\newcommand{\cL}{\ensuremath{\mathcal{L}}}

\newcommand{\cN}{\ensuremath{\mathcal{N}}}

\theoremstyle{plain}
\newtheorem{Thm}{Theorem}[section]

\newtheorem{Lem}[Thm]{Lemma}
\newtheorem{Prop}[Thm]{Proposition}

\theoremstyle{definition}
\newtheorem{Ass}[Thm]{Assumption}

\newtheorem{Rem}[Thm]{Remark}

\newcommand{\Proof}[2][Proof]{\begin{proof}[{#1}] #2 \end{proof}}

\numberwithin{equation}{section}

\makeatletter
\renewcommand\section{\@startsection {section}{1}{\z@}%
                                   {-3.5ex \@plus -1ex \@minus -.2ex}%
                                   {2.3ex \@plus.2ex}%
                                   {\normalfont\large\bf}}
\makeatother

\makeatletter
\renewcommand\subsection{\@startsection {subsection}{1}{\z@}%
                                   {-3.5ex \@plus -1ex \@minus -.2ex}%
                                   {2.3ex \@plus.2ex}%
                                   {\normalfont\normalsize\bf}}
\makeatother

\usepackage{color}
\renewcommand{\phi}{\hat{\psi}}

\usepackage{mathtools}
\mathtoolsset{showonlyrefs=true}

\allowdisplaybreaks[1]

\begin{document}

\begin{center}
{\Large \bf 
On the optimality of double barrier strategies for L\'evy processes
}
\end{center}
\begin{center}
Kei Noba
\end{center}

\begin{abstract}
This paper studies de Finetti's optimal dividend problem with capital injection. 
We confirm the optimality of a double barrier strategy when the underlying risk model follows a L\'evy process that may have positive and negative jumps. 
The main result in this paper is a generalization of \cite[Theorem 3]{AvrPalPis2007}, which is the spectrally negative case, and \cite[Theorem 3.1]{BayKypYam2013}, which is the spectrally positive case. 
In contrast with the spectrally one-sided cases, double barrier strategies cannot be handled by using scale functions to obtain some properties of the expected net present values (NPVs) of dividends and capital injections. 
Instead, to obtain these properties, we observe changes in the sample path (and the associated NPV) when there is a slight change to the initial value or the barrier value.
\end{abstract}
\section{Introduction}
In this paper, we study the optimal dividend problem with capital injection. 
A L\'evy process $X$ describes the risk process of an insurance company. 
The company pays dividends from the risk process $X$ and capital is injected into $X$ to avoid ruin. 
This paper aims to show a joint strategy for dividend payout and capital injection that maximizes the total expected dividend payments minus the cost of capital injections. 
\par
Here, we focus on the double barrier strategies. 
The double barrier strategy at $a\geq 0$ is the strategy in which (a) when the risk process exceeds the barrier $a$, 
the company pays a dividend determined by the excess over $a$, 
and (b) when the risk process falls below $0$, capital is injected accordingly to avoid ruin. 
The controlled process that results from applying a double barrier strategy behaves as a doubly reflected L\'evy process. 
\par
In previous studies of de Finetti's optimal dividend problem for spectrally negative L\'evy processes, 
the optimality of double barrier strategies has been proven. 
For example, Avram et al.\cite{AvrPalPis2007} proved optimality for general spectrally negative L\'evy processes. 
Furthermore, in expanded situations {that deal with two-sided singular control problems}, 
Baurdoux--Yamazaki\cite{BauYam2015} and 
Yamazaki\cite{Yam2018} proved optimality for general spectrally negative L\'evy processes. 
As with the spectrally negative case, 
some previous studies have considered spectrally positive L\'evy processes. 
Avanzi et al.\cite{AvaSheBer2011} proved the optimality of a double barrier strategy 
for spectrally positive compound Poisson processes. 
Bayraktar et al.\cite{BayKypYam2013} generalized the result to general spectrally positive L\'evy processes. 
In addition, many papers have considered other optimal dividend problems with capital injections for 
spectrally one-sided L\'evy processes 
(\cite{ZhaWanYaoChe2015}, 
\cite{PerYam2017_1}, \cite{ZhaWanYin2017}, \cite{ZhaCheYan2017}, \cite{PerYam2017_2}, 
\cite{PerYamYu2018}, \cite{NobPerYamYan2018}, \cite{CzaPerYam2018}, \cite{LiYin2019}). 
\par
Recently, de Finetti's optimal dividend problem for L\'evy processes with two-sided jumps has been studied. 
Some previous studies have considered de Finetti's optimal dividend problem without bail-outs. 
In particular, double or mixed-exponential jump diffusion processes have been discussed. 
For example, Bo et al.\cite{BoSonTanWanYan2012} computed the expected net present values (NPVs) of dividends of barrier strategies and gave 
numerical results for double exponential jump diffusion processes. 
Yin et al.\cite{YinSheWen2013} computed the expected NPVs of dividends of barrier strategies for mixed-exponential jump diffusion processes. 
Yuen--Yin\cite{YueYin2011} and Yin et al.\cite{YinYueShe2015} claim to have proven the optimality of the barrier strategy 
for more general L\'evy processes with two-sided jumps, but their proofs seem to have some flaws. 
In addition to these studies, Li et al.\cite{LiTanWanYan2014} gave some computational results 
that seem to provide the expected NPVs of dividends and capital injections of double barrier strategies for double-exponential jump diffusion processes. 
Overall, though, no existing paper seems to prove the optimality of any strategy.
\par
The objective of this paper is to show the optimality of a double barrier strategy for L\'evy processes that may have two-sided jumps. 
The class of L\'evy processes that we consider contains L\'evy processes with bounded variation paths and positive drifts, 
mixed-exponential jump diffusion processes, 
spectrally one-sided L\'evy processes, and others. 
\par
Our proof the optimality of a double barrier strategy has two broad steps. 
\begin{enumerate}
\item We select a candidate barrier $a^\ast \geq 0$ for double barrier strategies. 
In this step, we compute the derivative of $v_{\pi^a}(x)$ for $a$, where $v_{\pi^a}(x)$ is the expected NPV of dividends and capital injections of the double barrier strategy at $a\geq 0$ when the 
risk process $X$ starts from $x \in \bR$ (see Sections \ref{Sec04} and \ref{SecA}). 
\item To prove optimality for the chosen case, we apply a verification lemma to $v_{\pi^{a^\ast}}(x)$ as done in \cite{AvrPalPis2007}. 
Here, we need to find and use some properties of the derivative of $v_{\pi^a}(x)$ with respect to $x$
(see Sections \ref{Sec05}, \ref{Sec0B} and \ref{SecB}). 
\end{enumerate} 
A difficulty with this approach is how to obtain some properties of the derivative of the expected NPV of dividends and capital injections of double barrier strategies. 
In the case of spectrally one-sided L\'evy processes, we can represent the expected NPVs of dividends and capital injections of double barrier strategies by using scale functions, 
as done for \cite[(5.4)]{AvrPalPis2007} and \cite[(3.1)]{BayKypYam2013}. 
Since we know many properties of scale functions (see, e.g., \cite{KuzKypRiv2012} or \cite[Section 8]{Kyp2014}), 
we can obtain some properties of the derivative of each expected NPV. 
On the other hand, L\'evy processes that have two-sided jumps do not have scale functions associated with them. 
In the case of mixed-exponential jump diffusion processes, I predict we can represent the expected NPVs 
as the sum of exponential functions 
in the same way as done in \cite{BoSonTanWanYan2012} and \cite{YinSheWen2013}. 
However, I expect that the forms of the expected NPVs found in this way will be complicated to analyze. 
In addition, the expected NPVs cannot be expressed using this approach for more general Levy processes. 
So, we need to consider a new way to obtain some properties of the derivatives of the expected NPVs. 
In this paper, we obtain the derivatives of the expected NPVs 
by observing how the behavior of the sample path changes when either the initial value of the sample path or the value of the barrier is slightly shifted. 
Specifically, we represent the derivatives of the expected NPVs 
using the Laplace transforms of hitting times. 
\par
This paper is organized as follows. 
In Section \ref{Sec02}, we describe the notation and give some assumptions about L\'evy processes. 
In addition, we give the setting of the optimal dividend problem. 
In Section \ref{Sec03}, we give an overview of the double barrier strategies and confirm that they are admissible. 
In Section \ref{Sec04}, we select the candidate barrier $a^\ast$. 
In Section \ref{Sec05}, we prove the optimality of the double barrier strategy for $a^\ast$, using the verification lemma. 
The main result is in this section. 
In Section \ref{Sec06}, we give examples of L\'evy processes {with unbounded variation paths} that satisfy the assumptions given in Section \ref{Sec02}. 
In Section \ref{SecA}, we consider the behavior of doubly reflected L\'evy processes to compute the derivative of $v_{\pi^a}(x)$ with respect to $a$ and select the candidate barrier $a^\ast$. 
In Section \ref{Sec0B}, we give the proof of the verification lemma. 
In Section \ref{SecB}, we consider the behavior of doubly reflected L\'evy processes to compute the derivative of $v_{\pi^a}(x)$ with respect to $x$.

\section{Preliminalies}\label{Sec02}
\subsection{L\'evy processes}
In this section, we describe our notation and give some assumptions about the L\'evy processes considered in this paper. 
\par
Let $X=\{X_t : t \geq 0 \}$ be a L\'evy process defined on a probability space $(\Omega , \bF , \bP)$. 
For $x\in\bR$, we denote by $\bP_x$ the law of $X$ when it starts at $x$. 
Let $\Psi$ be the characteristic exponent of $X$ that satisfies 
\begin{align}
e^{-t\Psi(\lambda)}=\bE_0\sbra{e^{i\lambda X_t}},  \quad \lambda \in \bR, ~t\geq 0. 
\end{align}   
The characteristic exponent $\Psi$ is known to take the form  
\begin{align}
\Psi (\lambda) = -i\gamma\lambda +\frac{1}{2}\sigma^2 \lambda^2 
+\int_\bR (1-e^{i\lambda x}+i\lambda x1_{\{\absol{x}<1\}}) \Pi(dx) , ~~~~~~\lambda\in\bR. \label{202a}
\end{align}
Here, $\gamma\in\bR$, $\sigma\geq 0$, and $\Pi $ is a L\'evy measure on $\bR \backslash \{0\}$ such that 
\begin{align}
\int_{\bR\backslash \{ 0\}}(1\land x^2) \Pi (dx). 
\end{align}
The process $X$ has bounded variation paths if and only if $\sigma= 0$ and $\int_{\absol{x}<1} \absol{x}\Pi(dx)<\infty$. 
When this holds, we can write
\begin{align}
\Psi(\lambda)= -i\delta+\int_\bR(1-e^{i\lambda x})\Pi (dx),  \label{204}
\end{align}
where
\begin{align}
\delta = \gamma-\int_{\absol{x}<1}x \Pi(dx).
\end{align}
\par
Let $\cF=\{\cF_t : t\geq 0\}$ be the filtration generated by $X$. 
For $x\in \bR$, we write 
\begin{align}
\tau^+_x=\inf\{t> 0 : X_t > x\} \quad\text{and}\quad\tau^-_x=\inf\{ t>0: X_t < x\}.
\end{align}
We fix the discount factor $q>0$. 
For $a>0$ and $x \in \bR$, we write 
\begin{align}
\bar{\varphi}_{a, 0}= \bE_x \sbra{e^{-q\tau^+_a}; \tau^+_a < \tau^-_0}, \quad \un{\varphi}_{0, a} (x)=\bE_x \sbra{e^{-q\tau^-_0};   \tau^-_0 < \tau^+_a}. 
\end{align}
\par
For $a\in \bR$, let $Y^a$ be a reflected process defined by 
\begin{align}
Y^a_t = X_t - \rbra{ (\sup_{s\in [0, t]}X_s - a)  \lor 0 }, ~~~~~~t \geq 0. 
\end{align}
For $x \in \bR$, we write
\begin{align}
\kappa^{a, -}_x=\inf \{t>0 : Y^a_t < x\}. 
\end{align}
\par
We impose the following assumptions on $X$. 
\begin{Ass}\label{Ass101}
We assume that $X$ does not have monotone paths, 
and $X$ satisfies 
\begin{align}
\bE_0\sbra{\absol{X_1}}< \infty. \label{201a}
\end{align}
By \cite[Theorem 3.8]{Kyp2014}, the condition \eqref{201a} holds if and only if 
\begin{align}
\int_{\absol{x}\geq 1} \absol{x}  \Pi (dx) < \infty. \label{210}
\end{align}
If the process $X$ has unbounded variation paths, then we assume the maps $\bar{\varphi}_{a, 0}$ and $\un{\varphi}_{0, a}$ 
{have Radon--Nikodym densities $\bar{\varphi}_{a, 0}^\prime$ and $\un{\varphi}_{0, a}^\prime$ with respect to the Lebesgue measure,  which is continuous almost everywhere and locally bounded on $(0, a)$. }
\end{Ass}
\begin{Rem}
{Because the maps $\bar{\varphi}_{a, 0}$ and $\un{\varphi}_{0, a}$ are monotone functions, these are continuous on $[0, a]$ almost everywhere 
with respect to the Lebesgue measure. }
If the process $X$ has unbounded variation paths, it is easy to check that the maps $\bar{\varphi}_{a, 0}$ and $\un{\varphi}_{0, a}$ are continuous on $[0, a]$. 
\end{Rem}
\begin{Rem}
{
In fact,  the maps $\bar{\varphi}_{a, 0}$ and $\un{\varphi}_{0, a}$ are continuous on $[0, a]$ 
if $X$ is not a compound Poisson process. 
However, we do not give the proof of the fact since the fact is not important in this paper. 
}
\end{Rem}
{For $a>0$, we define 
\begin{align}
E^{(1)}_a= \{x \in [0, a]: \bar{\varphi}_{ a,0}\text{ or }\un{\varphi}_{0, a}\text{ are discontinuous at } x \}. 
\end{align} 
Note that $E^{(1)}_a$ is the null set. In addition, 
$E^{(1)}_a = \emptyset$ when $X$ has unbounded variation paths. }
\par
{
When $X$ has unbounded variation paths, we define 
\begin{align}
E^{(2)}_a= \{x \in (0, a): \bar{\varphi}_{ a,0}^\prime\text{ or }\un{\varphi}_{0, a}^\prime\text{ are discontinuous at } x \}.
\end{align}
Note that $E^{(2)}_a$ is the null set. 
}
\par
We define sets of functions $C^{(1)}_{\text{line}} $ and $C^{(2)}_{\text{line}} $. 
Let $C^{(1)}_{\text{line}}$ be the set of function 
$f\in C(\bR)$ {such that $f(y)-f(x)= \int_x^y f^\prime (u)du$ on $(0, \infty)$ 
for some locally bounded function $f^\prime$ on $(0, \infty)$ 
and such that $f$ satisfies} 
\begin{align}
\absol{f (x)} < b_1 \absol{x}+b_2, ~~~~~~x\in\bR,  \label{211b}
\end{align}
for some $b_1, b_2>0$. 
Let $C^{(2)}_{\text{line}}$ be the set of function 
$f$ in $C^{(1)}_{\text{line}}$ such that {$f$ is continuously differentiable on $(0, \infty)$, and }
$f^\prime(y)-f^\prime (x)=\int_x^y f^{\prime\prime}(u)du$ on $  (0, \infty)$ for some locally bounded function $f^{\prime\prime}$ on $(0 , \infty)$. 
Let $\cL$ be the operator applied to $f \in C^{(1)}_{\text{line}}$ (resp., $C^{(2)}_{\text{line}}$) for the case in which $X$ is of bounded (resp., unbounded) variation with 
\begin{align}
\cL f(x) = \gamma f^\prime (x)+\frac{1}{2}\sigma^2f^{\prime \prime}(x)+
\int_{\bR\backslash \{0\}}(f(x+z)-f(x)-f^\prime(x)z1_{\{\absol{z}<1\}})\Pi(dz),~x\in (0, \infty), \label{212b}
\end{align}
{for a fixed $f^\prime$ (resp., $f^{\prime\prime}$). }
\begin{Rem}\label{Rem203}
The integral in \eqref{212b} is well defined. We prove this fact here. 
We have 
\begin{align}
&\int_{\bR\backslash \{0\}}  \absol{f(x+z)-f(x)-f^\prime(x)z1_{\{\absol{z}<1\}}}\Pi(dz)\\
&\leq\int_{(-\infty , -1 \lor (\frac{ -x}{2})]\cup [1, \infty)}   \absol{f(x+z)-f(x)-f^\prime(x)z1_{\{\absol{z}<1\}}}\Pi(dz)\label{213}\\
&+\int_{(-1\lor (\frac{ -x}{2} ) , 1) \backslash \{0\}} \absol{f(x+z)-f(x)-f^\prime(x)z}\Pi(dz).\label{214}
\end{align}
By \eqref{210}, \eqref{211b}, and the definition of $\Pi$, \eqref{213} is finite. 
We have 
\begin{align}
\eqref{214}=&\int_{(-1\lor (\frac{ -x}{2}) , 1) \backslash \{0\}}\Pi(dz) \absol{ \int_0^z \rbra{f^\prime (x+y)- f^\prime(x)}dy}.  \label{216}
\end{align}
Because {$f^\prime$ is locally bounded, we have}
$\absol{f^\prime(x+ \cdot)} \leq b_3$ on $[-1\lor (\frac{ -x}{2}) , 1 ]$ for some $b_3 > 0$, we have 
\begin{align}
\eqref{216}\leq 2b_3 \int_{(-1\lor (\frac{ -x}{2}) , 1) \backslash \{0\}}\absol{z}\Pi(dz), 
\end{align}
and so \eqref{216} is finite when $X$ has bounded variation paths. 
When $X$ has unbounded variation paths, because $f^{\prime\prime}$ is locally bounded, we have 
\begin{align}
\absol{f^{\prime \prime}(y)}< b_4, \quad y \in (x+ (-1 \lor (\frac{ -x}{2})) , x+1) 
\end{align}
for some $b_4>0$. So, we have 
\begin{align}
\eqref{216}=&\int_{(-1\lor (\frac{ -x}{2}) , 1) \backslash \{0\}}\Pi(dz) \absol{ \int_0^z dy \int_0^y f^{\prime\prime} (x+w)dw}\\
\leq&\int_{(-1\lor (\frac{ -x}{2}) , 1) \backslash \{0\}}\Pi(dz)  \int_0^{\absol{z}} dy \int_0^y b_4dw\\
=&\frac{b_4}{2}\int_{(-1\lor (\frac{ -x}{2}) , 1) \backslash \{0\}} {\absol{z}}^2  \Pi(dz)< \infty. 
\end{align}
The proof is now complete. 
\end{Rem}
\begin{Rem}\label{Rem204}
By the proof in Remark \ref{Rem203}, and the dominated convergence theorem, it is easy to verify that the map 
\begin{align}
x &\mapsto \int_{\bR\backslash \{0\}}(f(x+z)-f(x))\Pi(dz) \\
\bigg{(}\text{resp., }x &\mapsto \int_{\bR\backslash \{0\}}(f(x+z)-f(x)-f^\prime(x)z1_{\{\absol{z}<1\}})\Pi(dz) \bigg{)}
\end{align}
is continuous on $(0, \infty)$ {when $X$ has bounded (resp., unbounded) variation paths}. 
\end{Rem}

\subsection{The optimal dividend problem with capital injection}
In this paper, a strategy is a pair of processes $\pi=\{(L^\pi_t, R^\pi_t): t\geq 0\}$ consisting of the cumulative amount of dividends $L^\pi$ and 
the cumulative amount of capital injection $R^\pi$. 
The corresponding risk process is given by $U^\pi_{0-}=X_0$, and 
\begin{align}
U^\pi_t=X_t-L^\pi_t +R^\pi_t,~~~~t\geq 0. 
\end{align}
\par
For the dividend strategy, we assume that $L^\pi$ is a non-decreasing, right-continuous, and $\cF$-adapted process with $L^{\pi}_{0-}=0$. 
\par
For the capital injection strategy, we assume that $R^\pi$ is a non-decreasing, right-continuous, and $\cF$-adapted process with $R^{\pi}_{0-}=0$, and 
\begin{align}
R^\pi_t \geq -(X_t-L^\pi_t), 
\quad t \geq 0. \label{230}
\end{align} 
The condition \eqref{230} implies that $U^{\pi}$ never hits $(-\infty , 0)$. 
\par
For $x\in\bR$, we write
\begin{align}
v^L_\pi (x)=\bE_x\sbra{\int_{[0, \infty)} e^{-qt}dL^\pi_t}, ~~
v^R_\pi (x)=\bE_x\sbra{\int_{[0, \infty)} e^{-qt}dR^\pi_t}.
\end{align}
Let $\beta>1$ be the cost per unit of injected capital. The objective is to maximize the expected NPV
\begin{align}
v_\pi (x)=v^L_\pi (x)-\beta v^R_\pi (x),~~~~x\in \bR 
\end{align}
over the set of all admissible strategies $\cA$ that satisfy all the constraints described above as well as $v^R_\pi (x)<\infty$ for $x\in\bR$. 
Hence, the problem is to obtain an optimal strategy $\pi^\ast$ satisfying 
\begin{align}
 v(x):= \sup_{\pi \in \cA} v_{\pi}(x)=v_{\pi^\ast} (x),~~~~x\in \bR. \label{108}
\end{align}

\section{Double barrier strategies}\label{Sec03}

The objective of this paper is to show the optimality of the double barrier strategy. 
In this section, we recall details of double barrier strategies, which are constructed 
in \cite[Section 4]{AvrPalPis2007}. \par
Let $a > 0$. The double barrier strategy $\pi^{a}$ at $a$ is the strategy constructed as follows. 
\begin{itemize}
\item[Step 0]
Set $T_0=\tau^-_0$, $T_a=\tau^+_a$, and $\eta =X_{T_0\land T_a}$. 
For $t<T_0\land T_a$, set $L^{\pi^a}_t=R^{\pi^a}_t=0$ and $U^{\pi^a}_t=X_t$.
If $T_0< T_a$, go to Step 2. If $T_0>T_a$, go to Step 1. 
\item[Step 1]
For $t\geq T_a$, we set 
\begin{align}
L^\prime_t = \eta-a+((X_t-X_{T_a})\lor 0),~~~~U^\prime_t =a+X_t-X_{T_a}-((X_t-X_{T_a})\lor 0). 
\end{align}
We reset $T_0= \inf\{ t> 0:U^\prime_t< 0 \}$ and $\eta=U^\prime_{T_0}$. 
For $t\in [T_a, T_0) $, we set $U^{\pi^a}_t =U^\prime_t$, 
$L^{\pi^a}_t =L^{\pi^a}_{T_a-}+L^\prime_t$, and $R^{\pi^a}_t =R^{\pi^a}_{T_a-}$. 
Go to Step 2. 
\item[Step 2]
For $t\geq T_0$, we set 
\begin{align}
R^\prime_t = -\eta-((X_t-X_{T_0})\land 0),~~~~U^\prime_t =X_t-X_{T_0}-((X_t-X_{T_0})\land 0). 
\end{align}
We reset $T_a= \inf\{ t> 0:U^\prime_t> a \}$ and $\eta=U^\prime_{T_a}$. 
For $t\in [T_0, T_a) $, we set $U^{\pi^a}_t =U^\prime_t$, $L^{\pi^a}_t =L^{\pi^a}_{T_0-}$, and $R^{\pi^a}_t =R^{\pi^a}_{T_0-}+R^\prime_t$. 
Go to Step 1. 
\end{itemize}
When $X$ has bounded variation paths, we can construct the double barrier strategy $\pi^0$ at $0$ as follows. For $t\geq 0$,  
\begin{align}
L^{\pi^0}_t&=X_01_{\{X_0>0\}} +\delta  t1_{\{\delta>0\}}+ \sum_{t\geq 0}(X_t-X_{t-})1_{\{ X_t-X_{t-}>0 \}}, \\
R^{\pi^0}_t&=-X_01_{\{X_0<0\}} -\delta  t1_{\{\delta<0\}}- \sum_{t\geq 0}(X_t-X_{t-})1_{\{ X_t-X_{t-}<0 \}}, \\
U^{\pi^0}_t&=0. 
\end{align}
\begin{Rem}
The doubly reflected processes are standard processes. 
\end{Rem}
\par
We prove the admissibility of the double barrier strategies in the following lemmas: 
\begin{Lem}\label{Lem201}
We have, for $a>0$ and $x \in \bR$, 
\begin{align}
 v^L_{\pi^a}(x) < \infty, \quad v^R_{\pi^a}(x) < \infty. 
\end{align}
\end{Lem}
\Proof{
The proof of $v^L_{\pi^a}(x)<\infty$ is the same as the proof of $v^R_{\pi^a}(x)<\infty$, so we prove only $v^R_{\pi^a}(x)<\infty$. 
\par
By the definition of $\pi^a$, 
we have 
\begin{align}
v_{\pi^a}^R(x) =
\begin{cases}
-x+v_{\pi^a}^R (0), \quad &x<0, \\
v_{\pi^a}^R (a), \quad & x>a.
\end{cases} 
\end{align}
Thus it is sufficient to prove that 
\begin{align}
\bar{v}_{\pi^a}^R:= \sup_{x\in [0, a]} v_{\pi^a}^R (x) < \infty . \label{203a}
\end{align}
\par
We write $R_t=-\rbra{(\inf_{s\in[0, t]}X_s)\land 0 }$  
and fix $u>0$. 
For $n\in\bN$, we define 
\begin{align}
\mu^{[n]}_a= (u+\mu^{[n-1]}_a)\land \inf \{t>\mu^{[n-1]}_a : U^{\pi^a}_t=a, \text{ there exists }s \in [\mu^{[n-1]}_a, t) \text{ such that }U_s = 0 \}, 
\end{align}
where $\mu^0_a =0$. 
By the strong Markov property, we have 
\begin{align}
\bar{v}_{\pi^a}^R
&= \lim_{n\uparrow \infty }
\sup_{x \in [0,a]} \bE_x \sbra{\int_0^{\mu^{[n]}_a}e^{-qt} dR^{\pi^a}_t} \\
&= \lim_{n\uparrow \infty }\sup_{x\in[0, a]}\rbra{\sum_{k=1}^n\bE_x \sbra{\int_{\mu^{[k-1]}_a}^{\mu^{[k]}_a}e^{-qt} dR^{\pi^a}_t} }\\
&\leq\lim_{n\uparrow \infty }\rbra{    \sum_{k=1}^n\rbra{ \bE_0 \sbra{ e^{-q\mu^{[1]}_a}}}^{k-1}   \sup_{x\in[0,a]}
  \bE_x \sbra{\int_{0}^{\mu^{[1]}_a}e^{-qt} dR^{\pi^a}_t}    }\\
  &=\frac{ \sup_{x \in[0,a]}\bE_x \sbra{\int_{0}^{\mu^{[1]}_a}e^{-qt} dR^{\pi^a}_t} }{1- \bE_0\sbra{e^{-q\mu^{[1]}_a}}},  \label{313}
\end{align}
and by the definition of $\mu^{[1]}_a$, we have 
 \begin{align}
 \eqref{313}
\leq \frac{\sup_{x \in[0,a]}\bE_x\sbra{ \int_0^{u} e^{-qt} dR_t }}{1- \bE_0 \sbra{ e^{-q \mu^{[1]}_a } }} 
 \leq \frac{\bE_0\sbra{ \absol{\inf_{t \in [0, u] }X_t }}}{1- \bE_0 \sbra{ e^{-q \mu^{[1]}_a } }}.  \label{207b}
\end{align}
\par
By the L\'evy--It\^o decomposition, there exists some $\delta^X \in \bR$, a zero-mean square-integrable martingale $M^X$ starting from $0$, a Poisson process $N^X$, and 
a sequence of independent and identically distributed (i.i.d.) random variables ${\{ J^X_n\}}_{n\in \bN}$ taking values in $(-\infty, -1) \cup (1, \infty)$ such that 
\begin{align}
X_t = \delta^X t + M^X_t+ \sum_{i=1}^{N^X_t} J^X_i,  \quad t\geq 0, \quad \bP_0\text{-a.s.}         \label{208}
\end{align}
By Doob's 
maximal inequality, we have 
\begin{align}
\bE_0 \sbra{ \absol{ \inf_{t \in [0, u] } M^X_t  }} \leq 
1+ \bE_0 \sbra{{ \rbra{ \inf_{t \in [0, u] } M^X_t }}^2}\leq 1+ 4\bE_0 \sbra{{(M^{X}_u)}^2}<\infty . \label{209a}
\end{align}
By the compensation theorem of the Poisson point processes and \eqref{210}, we have 
\begin{align}
\bE_0 \sbra{   \absol{  \inf_{t\in [0, u]}\sum_{i=1}^{N^X_t} J^X_i } }
&\leq \bE_0 \sbra{   \absol{  \sum_{i=1}^{N^X_u} (J^X_i \land 0)} }\\
&=-\bE_0\sbra{\int_{[0, \infty) \times (-\infty, -1)}e^{-qt} x\cN(dt \times dx) } \\
&=-\int_0^\infty e^{-qt}dt \int_{(-\infty,-1)}x \Pi (dx) < \infty, \label{211a}
\end{align}
where $\cN$ is a Poisson random measure 
on $([0, \infty) \times\bR ,\cB [0, \infty ) \times \cB (\bR), ds \times \Pi(dx)) $ 
associated with the jumps of $X$. 
By \eqref{208}, \eqref{209a}, and \eqref{211a}, we have 
\begin{align}
\bE_0 \sbra{   \absol{  \inf_{t\in [0, u]}X_t } }
\leq & \absol{\delta^X}u+ \bE_0\sbra{ \absol{  \inf_{t\in [0, u]}M^X_t } }+\bE_0 \sbra{   \absol{  \inf_{t\in [0, u]}\sum_{i=1}^{N^X_t} J^X_i } }< \infty. \label{213b}
\end{align}
By \eqref{207b} and \eqref{213b}, we obtain \eqref{203a}, and the proof is complete. 
}
\begin{Lem}
If $X$ has bounded variation paths, then for $x \in \bR$ we have 
\begin{align}
 v^L_{\pi^0}(x) < \infty, \quad  v^R_{\pi^0}(x) < \infty. 
\end{align}
\end{Lem}
\Proof{
By the same argument as that in the proof of Lemma \ref{Lem201}, it is sufficient to prove that $v^R_{\pi^0}(0)<\infty$. 
By the definition of $\pi^0$,  we have 
\begin{align}
v^R_{\pi^0}(0)&=-1_{\{\delta<0\}}\delta\int_0^\infty e^{-qt}dt-\bE_0\sbra{\sum_{t>0}e^{-qt}(X_t-X_{t-})1_{\{X_t-X_{t-}<0\}}  }.
\end{align}
By the compensation theorem of the Poisson point processes, we have 
\begin{align}
-\bE_0\sbra{\sum_{t>0}e^{-qt}(X_t-X_{t-})1_{\{X_t-X_{t-}<0\}}  }
&=-\bE_0\sbra{\int_{[0, \infty) \times (-\infty, 0)}e^{-qt} x\cN(dt \times dx) } \\
&=-\int_0^\infty e^{-qt}dt \int_{(-\infty,0)}x \Pi (dx). \label{325}
\end{align}
Because $X$ has bounded variation paths and \eqref{210} holds, we have $\eqref{325}< \infty$. 
The proof is complete. 
}
\begin{Rem}
We note that the assumption \eqref{201a} is necessary to prove the optimality 
of a double barrier strategy. We show this briefly here. 
Suppose $X$ has bounded (resp., unbounded) variation paths. 
If $\bE_0\sbra{X_1 \land 0}=-\infty$, then we can prove that $v^R_{\pi^a}(x)=\infty$ for $a\in [0, \infty)$ (resp., $a\in (0,\infty)$) and $x\in\bR$, and 
so no double barrier strategies are admissible. 
On the other hand, if $\bE_0\sbra{X_1 \lor 0}=\infty$, then we can prove that $v^L_{\pi^a}(x)=\infty$ for $a\in [0, \infty)$ and $x\in\bR$, 
and so this problem does not make sense. 
\end{Rem}

\section{Selection of the candidate barrier $a^\ast$}\label{Sec04}
In this section, we focus on the double barrier strategy defined in the previous section and choose a candidate barrier $a^\ast$. 
\par
We define 
\begin{align}
a^\ast=\inf \cbra{a>0: \beta \un{\nu}(a) \leq1}, 
\end{align}
where $\un{\nu}(a)=\bE_a \sbra{ e^{-q\kappa^{a, -}_{0}}}$. 
Since $X$ has stationary independent increments, the map $a\mapsto \un{\nu}(a)$ is non-increasing, and we have $\lim_{a \uparrow \infty }\un{\nu}(a)=0 $, 
so $a^\ast < \infty$.  
\begin{Lem}\label{Lem202}
For $a, x \in [0, \infty)$, we have 
\begin{align}
v_{\pi^{a^\ast}}(x ) \geq v_{\pi^a}(x). 
\end{align}
\end{Lem}
We postpone the proof of Lemma \ref{Lem202} until after the proof of Lemma \ref{Lem302a}. 
\par
For $n\in\bN$, we define the hitting times for $n \in \bN$ by induction as follows. For $n\in\bN$, 
\begin{align}
\bar{\rho}^{a, [n]}_{a}=\inf\{t\geq \un{\rho}^{a, [n-1]}_{0} : L^{\pi^a}_t >L^{\pi^a}_{\un{\rho}^{a, [n-1]}_{0}-}\}, \quad
\un{\rho}^{a, [n]}_{0}=\inf\{t\geq \bar{\rho}^{a, [n]}_{a} : R^{\pi^a}_t >R^{\pi^a}_{\bar{\rho}^{a, [n]}_{a}-}\}, 
\end{align}
where $\un{\rho}^{a, (0)}_{0}=0$.
For simplicity, we write $ \bar{\nu}_x(a)=\bE_x \sbra{ e^{-q\bar{\rho}^{a, [1]}_{a}}}$. 
To compute the derivative of $v_{\pi^a}(x)$ with respect to $a$, 
we write $V_x(a):=v_{\pi^a}(x)$, $V^L_x(a):=v^L_{\pi^a}(x)$ and $V^R_x(a):=v^R_{\pi^a}(x)$. 
\begin{Lem}\label{Lem302}
For $x \in \bR$ and $a \in (0, \infty)$, we have 
\begin{align}
&\lim_{\epsilon\downarrow 0} \frac{V_x^L(a+\epsilon)-V_x^L (a)}{\epsilon}
=\frac{-\bar{\nu}_x (a)}{1-\un{\nu} (a)\bar{\nu}_0 (a) }, \label{304}\\
&
\lim_{\epsilon\downarrow 0} \frac{V_x^R(a+\epsilon)-V_x^R (a)}{\epsilon}=
\frac{-\bar{\nu}_x (a)\un{\nu} (a)}{1-\un{\nu} (a)\bar{\nu}_0 (a) }. \label{305}
\end{align}
\end{Lem}
\Proof{
We  estimate $\bar{\lim}_{\epsilon\downarrow 0} \frac{V_x^L(a+\epsilon)-V_x^L (a)}{\epsilon}$ and $\un{\inf}_{\epsilon\downarrow 0} \frac{V_x^L(a+\epsilon)-V_x^L (a)}{\epsilon}$. 
We have
\begin{align}
V_x^L(a)-V_x^L (a+\epsilon)&=\bE_x\sbra{\int_{[0, \infty)} e^{-qt} d\rbra{ L^{\pi^a}_t -L^{\pi^{a+\epsilon}}_t  } }\\
&\begin{aligned}
=&\sum_{k=1}^\infty \bE_x\sbra{\int_{[\un{\rho}^{a, [k-1]}_0 , \bar{\rho}^{a, [k]}_a)} e^{-qt} d\rbra{ L^{\pi^a}_t -L^{\pi^{a+\epsilon}}_t  } }\\
&+\sum_{k=1}^\infty \bE_x\sbra{\int_{[\bar{\rho}^{a, [k]}_a , \un{\rho}^{a, [k]}_0)} e^{-qt} d\rbra{ L^{\pi^a}_t -L^{\pi^{a+\epsilon}}_t  } } \label{307}
\end{aligned}
\\
&\begin{aligned}
=&\sum_{k=1}^\infty \bE_x\sbra{\int_{[\un{\rho}^{a+\epsilon, [k-1]}_0 , \bar{\rho}^{a+\epsilon, [k]}_{a+\epsilon}]} e^{-qt} d\rbra{ L^{\pi^a}_t -L^{\pi^{a+\epsilon}}_t  } }\\
&+\sum_{k=1}^\infty \bE_x\sbra{\int_{(\bar{\rho}^{a+\epsilon, [k]}_{a+\epsilon} , \un{\rho}^{a+\epsilon, [k]}_0)} e^{-qt} d\rbra{ L^{\pi^a}_t -L^{\pi^{a+\epsilon}}_t  } }. \label{309}
\end{aligned}
\end{align}
The behaviors of $L^{\pi^a}$ and $L^{\pi^{a+\epsilon}}$ are given in Section \ref{SecA}. 
Because we have \eqref{414a} and \eqref{A07}, for $t \in [\un{\rho}^{a, [k-1]}_0 , \bar{\rho}^{a, [k]}_a)$, we have, for $k\in\bN$, 
\begin{align}
 \bE_x\sbra{\int_{[\un{\rho}^{a, [k-1]}_0 , \bar{\rho}^{a, [k]}_a)} e^{-qt} d\rbra{ L^{\pi^a}_t -L^{\pi^{a+\epsilon}}_t  } }=0. \label{310}
\end{align}
Because we have \eqref{427a} for $t \in[\bar{\rho}^{a, [k]}_a , \un{\rho}^{a, [k]}_0)$, we have, for $k\in\bN$,
\begin{align}
\bE_x\sbra{\int_{[\bar{\rho}^{a, [k]}_a , \un{\rho}^{a, [k]}_0)} e^{-qt} d\rbra{ L^{\pi^a}_t -L^{\pi^{a+\epsilon}}_t  } }\leq \epsilon \bE_x\sbra{e^{-q\bar{\rho}^{a, [k]}_a}}. \label{311}
\end{align}
By \eqref{307}, \eqref{310}, and \eqref{311}, we have 
\begin{align}
V_x^L(a)-V_x^L (a+\epsilon)\leq& \epsilon \sum_{k=1}^\infty\bE_x\sbra{e^{-q\bar{\rho}^{a, [k]}_a}} \\
=&\epsilon\bE_x\sbra{e^{-q\bar{\rho}^{a, [1]}_a}} \sum_{k=0}^\infty  \rbra{\bE_a\sbra{e^{-q\kappa^{a,-}_0 }}  \bE_0\sbra{e^{-q\bar{\rho}^{a, [1]}_a}}}^k\\
=&\epsilon\frac{\bar{\nu}_x (a)}{1-\un{\nu} (a)\bar{\nu}_0 (a) } .\label{314}
\end{align}
\par
By \eqref{418a} and \eqref{425a}, we have 
\begin{align}
0\leq U^{\pi^{a+\epsilon}}_t-U^{\pi^a}_t \leq \epsilon ,\quad t\geq 0. 
\end{align}
So, we have 
\begin{align}
U^{\pi^{a+\epsilon}}_{\un{\rho}^{a+\epsilon, [k-1]}_0}-U^{\pi^a}_{\un{\rho}^{a+\epsilon, [k-1]}_0}=0-0=0, 
\quad 
U^{\pi^{a+\epsilon}}_{\bar{\rho}^{a+\epsilon, [k]}_{a+\epsilon}}-U^{\pi^a}_{\bar{\rho}^{a+\epsilon, [k]}_{a+\epsilon}}=(a+\epsilon)-a=\epsilon.   
\end{align}
By \eqref{417a} and \eqref{A16}
, $R^{\pi^a}-R^{\pi^{a+\epsilon}}$ is non-decreasing. 
Therefore, 
$L^{\pi^a}-L^{\pi^{a+\epsilon}}$ increases by at least $\epsilon$ on $[\un{\rho}^{a+\epsilon, [k-1]}_0 ,\bar{\rho}^{a+\epsilon, [k]}_{a+\epsilon} ]$. 
By \eqref{A07} and \eqref{427a}, $L^{\pi^a}-L^{\pi^{a+\epsilon}}$ is non-decreasing, giving us 
\begin{align}
 \bE_x\sbra{\int_{[\un{\rho}^{a+\epsilon, [k-1]}_0 , \bar{\rho}^{a+\epsilon, [k]}_{a+\epsilon}]} e^{-qt} d\rbra{ L^{\pi^a}_t -L^{\pi^{a+\epsilon}}_t  } } 
 \geq\epsilon \bE_x\sbra{e^{-q\bar{\rho}^{a+\epsilon, [k]}_{a+\epsilon}}} . \label{315}
\end{align}
By \eqref{309} and \eqref{315}, we have 
\begin{align}
V_x^L(a)-V_x^L (a+\epsilon)\geq& \epsilon \sum_{k=1}^\infty\bE_x\sbra{e^{-q\bar{\rho}^{a+\epsilon, [k]}_{a+\epsilon}}} 
=\epsilon\frac{\bar{\nu}_x (a+\epsilon)}{1-\un{\nu} (a+\epsilon)\bar{\nu}_0 (a+\epsilon) } .\label{317}
\end{align}
It is easy to check that $\un{\nu}$ and $\bar{\nu}_x$ are right continuous. 
Thus, by \eqref{314} and \eqref{317}, we obtain \eqref{304}. 
\par
By the same argument as above, we obtain \eqref{305}. The proof is complete. 
}
\begin{Lem}\label{Lem302a}
Suppose $X$ has bounded variation paths. Then for $x \in \bR$ we have 
\begin{align}
\lim_{a\downarrow 0}V_x(a)=V_x (0). 
\end{align}
\end{Lem}
\Proof{We assume that the drift parameter $\delta$ is less than $0$. 
We prove $\lim_{a\downarrow 0}V^L_x(a)=V^L_x (0)$ and $\lim_{a\downarrow 0}V^R_x(a)=V^R_x (0)$. 
By the construction of $\pi^a$, $L^{\pi^a}$ increases only when $X$ takes positive jumps, and we have 
\begin{align}
(X_0-a)1_{\{X_0 -a >0\}} + \sum_{t\geq 0}e^{-qt}(X_t-X_{t-}-a)1_{\{ X_t-X_{t-}>a \}} \leq L^{\pi^a}_t \leq L^{\pi^0}_t, \quad t\geq 0. 
\end{align}
Because 
\begin{align}
\lim_{a\downarrow 0}\rbra{(X_0-a)1_{\{X_0 -a >0\}} + \sum_{t\geq 0}e^{-qt}(X_t-X_{t-}-a)1_{\{ X_t-X_{t-}>a \}}}=L^{\pi^0}_t, \quad t\geq 0 ,
\end{align}
we have 
\begin{align}
\lim_{a\downarrow 0} L^{\pi^a}_t =L^{\pi^0}_t , \quad t\geq 0, \quad \text{ and thus } \quad \lim_{a\downarrow 0}V^L_x(a)=V^L_x(0). \label{322}
\end{align}
By the construction of $\pi^a$, we can write 
\begin{align}
R^{\pi^a}_t =
-\rbra{\inf_{s\leq t }(X_s -L^{\pi^a}_s)} \land 0 ,\quad t\geq 0. 
\end{align}
By \eqref{322}, we have 
\begin{align}
\lim_{a\downarrow0}(X_t -L^{\pi^a}_t)=X_01_{\{X_0<0\}}+\delta t 1_{\{\delta < 0\}}+ \sum_{t\geq 0}(X_t-X_{t-})1_{\{ X_t-X_{t-}<0 \}}, 
\end{align}
so we have 
\begin{align}
\lim_{a\downarrow 0} R^{\pi^a}_t =R^{\pi^0}_t , \quad t\geq 0, \quad \text{ and thus } \quad \lim_{a\downarrow 0}V^R_x(a)=V^R_x(0).
\end{align}
\par
When $\delta$ is assumed to be positive, the proof is almost the same as above. 
The proof is now complete. 
}
\Proof[Proof of Lemma \ref{Lem202}]{
Since $\frac{-\bar{\nu}_x (a)}{1-\un{\nu} (a)\bar{\nu}_0 (a) }$ and $\frac{-\bar{\nu}_x (a)\un{\nu} (a)}{1-\un{\nu} (a)\bar{\nu}_0 (a) }$ are non-decreasing, and 
by Lemma \ref{Lem302}, both $V^L_x$ and $V^R_x$ are convex functions having Radon--Nikodym densities
\begin{align}
V_x^{L\prime}(a)=\frac{-\bar{\nu}_x (a)}{1-\un{\nu} (a)\bar{\nu}_0 (a) } ,\quad
V_x^{R\prime}(a)=\frac{-\bar{\nu}_x (a)}{1-\un{\nu} (a)\bar{\nu}_0 (a) }  \un{\nu} (a)
,\quad a\in (0, \infty)
\end{align}
with respect to the Lebesgue measure. 
This implies that $V_x $ has Radon--Nikodym density
\begin{align}
V_x^\prime(a)=\frac{-\bar{\nu}_x (a)}{1-\un{\nu} (a)\bar{\nu}_0 (a) } \rbra{1- \beta \un{\nu} (a)},\quad a\in (0, \infty). 
\end{align}
By the definition of $a^\ast$, $V_x$ is non-decreasing on $(0, a^\ast)$ and non-increasing on $(a^\ast , \infty)$. 
By Lemma \ref{Lem302a}, the proof is complete. 
}

\section{Verification}\label{Sec05}
In this section, we show the optimality of the strategy $\pi^{a^\ast}$ for the value $a^\ast$ selected in the previous section. 
The main theorem is the following. 
\begin{Thm}\label{Thm501}
The strategy $\pi^{a^\ast}$ is optimal, and the value function of the problem \eqref{108} is given by $v=v_{\pi^{a^\ast}}$.
\end{Thm}
We postpone the proof of Theorem \ref{Thm501} until after the proof of Lemma \ref{Lem408}. 
To prove Theorem \ref{Thm501}, we use the following verification lemma. 
\begin{Prop}
\label{Prop201}
Suppose that $X$ has bounded (resp., unbounded) variation paths. 
Let $w$ be a function on $\bR$ belonging to $ C^{(1)}_{\text{line}}$ (resp., $C^{(2)}_{\text{line}}$) and satisfying 
\begin{align}
\cL w (x)-q w (x) \leq 0,       ~~~~~~  &x > 0, \label{201}\\ 
1 \leq w^\prime (x) \leq \beta, ~~~~~~&x \in \bR \backslash \{0\}. \label{202}
\end{align}
Then we have $w(x)\geq v (x)$ for all $x\geq0$. 
\end{Prop}
We give the proof of Proposition \ref{Prop201} in Section \ref{Sec0B}. 
The purpose of this section is to prove that the strategy $\pi^{a^\ast}$ satisfies the conditions in Proposition \ref{Prop201}. 
\par
To apply Proposition \ref{Prop201} to $v_{\pi^\ast}$, we give Lemmas \ref{Thm203}, \ref{Lem404}, and \ref{Lem203}, 
which are lemmas about the smoothness of $v_{\pi^{a^\ast}}$ and some properties of the derivative of $v_{\pi^{a^\ast}}$.  
\begin{Lem}\label{Thm203}
We fix $a>0$. 
For $x\in(0, a)\backslash E^{(1)}_a$, we 
{can take a Radon--Nikodym density of $v_{\pi^{a}}$ as }  
\begin{align}
v^{\prime}_{\pi^{a}}(x)=\bar{\varphi}_{a, 0} (x) 
+\beta \un{\varphi}_{0, a}(x)
. \label{207}
\end{align}
\end{Lem}
\Proof{
We compute the derivative of $v^L_{\pi^a}$. 
For $y\in \bR$ and $t \geq 0$, we write $X^{(y)}_t=X_t+y$. 
For $b \in \bR$, we write $\tau^{(y),+}_b = \inf\{ t> 0: X^{(y)}_t \geq b\}$ and $\tau^{(y),-}_b = \inf\{ t> 0: X^{(y)}_t \leq b\}$. 
We write $L^{(y)}$ for the process that represents the cumulative amount of dividends of $X^{(y)}$ on which 
the double barrier strategy at $a$ has been imposed. 
We write $U^{(y)}$ for its surplus process. 
\par
{For $x \in [0, a-\epsilon]$, }
we can rewrite $v^L_{\pi^a}(x+\epsilon)-v^L_{\pi^a}(x)$ as 
\begin{align}
v^L_{\pi^a}(x+\epsilon)-&v^L_{\pi^a}(x)
=\bE_0\sbra{\int_{[0, \infty)} e^{-qt}d(L^{(x+\epsilon)}_t-L^{(x)}_t) }. \label{346}
\end{align}
The behaviors of $U^{(x)}$ and $U^{(x+\epsilon)}$ are summarized in Section \ref{SecB}. 
From Section \ref{SecB}, we have 
\begin{align}
&L^{(x+\epsilon)}_t-L^{(x)}_t \in [0, \epsilon] \text{ is non-decreasing for }t\geq 0; \label{506}\\
&R^{(x)}_t-R^{(x+\epsilon)}_t \in [0, \epsilon] \text{ is non-decreasing for }t\geq 0; \label{507}\\
&U^{(x+\epsilon)}_t-U^{(x)}_t=\epsilon- (L^{(x+\epsilon)}_t-L^{(x)}_t)-(R^{(x)}_t-R^{(x)}_t)\geq 0 \text{ for }t\geq 0. \label{508}
\end{align}
\par
Because $L^{(x+\epsilon)}_{\tau^{(x+\epsilon),+}_a-}=0$ and \eqref{506} holds, we have 
\begin{align}
\inf\{ t>0 : L^{(x+\epsilon)}_t -L^{(x)}_t >0 \}\geq \tau^{(x+\epsilon),+}_a. \label{407}
\end{align}
On $\{ \tau^{(x+\epsilon),+}_a >\tau^{(x+\epsilon),-}_0\}$, 
we have, by \eqref{407} and \eqref{508},
\begin{align}
L^{(x+\epsilon)}_{\tau^{(x+\epsilon),-}_0}-L^{(x)}_{\tau^{(x+\epsilon),-}_0}=0, \quad 
U^{(x+\epsilon)}_{\tau^{(x+\epsilon),-}_0}=U^{(x)}_{\tau^{(x+\epsilon),-}_0}=0, \quad 
R^{(x)}_{\tau^{(x+\epsilon),-}_0}-R^{(x+\epsilon)}_{\tau^{(x+\epsilon),-}_0}=\epsilon, 
\end{align}
and so by \eqref{507} and \eqref{508}, we have
\begin{align}
\inf\{ t>0 : L^{(x+\epsilon)}_t -L^{(x)}_t >0 \}=\infty. \label{408}
\end{align}
By \eqref{346}, \eqref{506}, \eqref{407}, and \eqref{408}, 
we have 
\begin{align}
v^L_{\pi^{a}}(x+\epsilon)-v^L_{\pi^{a}}(x)&\leq
\epsilon\bE_0 \sbra{e^{-q\tau^{(x+\epsilon),+}_a } ;\tau^{(x+\epsilon),+}_a <\tau^{(x+\epsilon),-}_0}
= \epsilon  \bar{\varphi}_{a, 0} (x+\epsilon). 
\label{410}
\end{align}
\par
Because $R^{(x)}_{\tau^{(x),-}_0-}=0$ and \eqref{507} holds, we have 
\begin{align}
\inf \{t>0 :R^{(x)}_t - R^{(x+\epsilon)}_t >0\}\geq \tau^{(x), -}_0. \label{514}
\end{align}
On $\{ \tau^{(x),+}_a <\tau^{(x),-}_0\}$, by \eqref{514} and \eqref{508}, we have 
\begin{align}
R^{(x)}_{\tau^{(x),+}_a} - R^{(x+\epsilon)}_{\tau^{(x),+}_a}=0, \quad U^{(x+\epsilon)}_{\tau^{(x),+}_a} = U^{(x)}_{\tau^{(x),+}_a}=a, \quad 
L^{(x+\epsilon)}_{\tau^{(x),+}_a} - L^{(x)}_{\tau^{(x),+}_a}= \epsilon, 
\end{align}
which implies that 
\begin{align}
 \inf\{ t>0 : L^{(x+\epsilon)}_t -L^{(x)}_t =\epsilon \}\leq \tau^{(x),+}_a . \label{411}
\end{align}
By \eqref{506} and \eqref{411}, 
we have 
\begin{align}
 v^L_{\pi^{a}}(x+\epsilon)-v^L_{\pi^{a}}(x)
&\geq \epsilon \bE_0 \sbra{e^{-q\tau^{(x),+}_a } ;\tau^{(x),+}_a <\tau^{(x),-}_0} 
=\epsilon   \bar{\varphi}_{a, 0} (x)
. \label{414}
\end{align}
Following the logic of the proofs of \eqref{410} and \eqref{414}, we have, {for $x \in [\epsilon , a]$,}
\begin{align}
\epsilon \bar{\varphi}_{a, 0} (x -\epsilon)
\leq v^L_{\pi^{a}}(x)-v^L_{\pi^{a}}(x-\epsilon)
\leq \epsilon  \bar{\varphi}_{a, 0}  (x). 
\label{415}
\end{align}
From \eqref{410}, \eqref{414}, and \eqref{415}, {for $x\in(0, a)\backslash E^{(1)}_a$, }we obtain 
\begin{align}
\lim_{\epsilon \downarrow 0 }\frac{v^L_{\pi^{a}}(x+\epsilon)-v^L_{\pi^{a}}(x)}{\epsilon}
=\lim_{\epsilon \downarrow 0 }\frac{v^L_{\pi^{a}}(x)-v^L_{\pi^{a}}(x-\epsilon)}{\epsilon}
=\bar{\varphi}_{a, 0} (x) .
\label{212a}
\end{align}
By a computation similar to that for the derivative of $v^L_{\pi^a}$, {for $x\in(0, a)\backslash E^{(1)}_a$, }
we obtain 
\begin{align}
\lim_{\epsilon \downarrow 0 }\frac{v^R_{\pi^{a}}(x+\epsilon)-v^R_{\pi^{a}}(x)}{\epsilon}
=\lim_{\epsilon \downarrow 0 }\frac{v^R_{\pi^{a}}(x)-v^R_{\pi^{a}}(x-\epsilon)}{\epsilon}=
-\un{\varphi}_{0, a} (x). 
\end{align}
{
In addition, $v^{L}_{\pi^{a}}$ and $v^{R}_{\pi^{a}}$ are continuous concave functions on $[0, a]$ from \eqref{410}, \eqref{414}, \eqref{415}, 
and the other computations, 
and so $v_{\pi^{a}}$ has a Radon--Nikodym density \eqref{207}. }
The proof is complete. 
}
\begin{Lem}\label{Lem404}
For $a\geq 0$, the function $v_{\pi^a}$ is a continuous function. 
\end{Lem}
\Proof{
By the definition of $\pi^a$, we have 
\begin{align}
v_{\pi^a}(x)=
\begin{cases}
v_{\pi^a}(a)+(x-a),\quad&x\geq a, \\
v_{\pi^a}(0)+\beta x,\quad&x\leq 0.
\end{cases}
\label{417}
\end{align}
So $v_{\pi^a}$ is continuous 
on $(-\infty, 0]\cup[a, \infty)$. 
{In addition, $v_{\pi^a}$ is continuous on $[0, a]$ by the proof of Lemma \ref{Thm203}. }
The proof is complete. 
}
{Since we have \eqref{417}, 
we define $v_{\pi^{a^\ast}}^\prime$ as 
\begin{align}
v_{\pi^{a^\ast}}^\prime(x)=&
\begin{cases}
1, \quad x \geq 1, \\
\beta, \quad x \leq 0.
\end{cases}
\label{520}\\
v_{\pi^a}^{\prime\prime}(x)=&0, \quad x \in (-\infty, 0]\cup[a, \infty). 
\label{544}
\end{align}
}
\begin{Lem}\label{Lem203}
We have 
\begin{align}
&1\leq v_{\pi^{a^\ast}}^\prime (x) \leq \beta, \quad  x \in (0, a^\ast)\backslash E^{(1)}_{a^\ast}, \label{419}
\end{align}
{and $v_{\pi^{a^\ast}}$ is a concave function on $(0, \infty)$. 
In addition, we have 
\begin{align}
v_{\pi^{a^\ast}}^\prime (a^\ast-)=1  \label{420}
\end{align}
when $X$ has unbounded variation paths. 
}
\end{Lem}
\Proof{
{\textbf{\textit{i})}}
It is easy to check that $\un{\nu}$ is right-continuous, and so we have $\beta\un{\nu}(a^\ast) \leq 1$. 
In this step, we define a constant $p^\ast$, and stopping times $K^{p^\ast}_0$ and $T^{p^\ast}_0$ 
for two cases. One is the case in which $\beta\un{\nu}(a^\ast) = 1$. The other is the case in which $\beta\un{\nu}(a^\ast) < 1$. 
\par
Suppose $\beta\un{\nu}(a^\ast) = 1$. Then we define $p^\ast=1$, $K^{p^\ast}_0=K^1_0= \kappa^{a^\ast , -}_0 $, and $T^{p^\ast}_0=T^1_0= \tau^-_0$. 
Here, we have 
\begin{align}
\beta\bE_{a^\ast}\sbra{ e^{-q K^{p^\ast}_0} }=1. \label{458}
\end{align}
\par
Suppose $\beta\un{\nu}(a^\ast) <1$. Then $X$ has bounded variation paths. 
For $n\in\bN$, we write 
\begin{align}
T[n]=&\inf\{t> T[n-1]: {Y^{a^\ast}_s\neq 0 \text{ for some }s\in (T[n-1], t),} ~Y^{a^\ast}_t=0 \}, \\
S[n]=&\inf\{t> S[n-1]: {X_s\neq 0 \text{ for some }s\in (S[n-1], t),}~X_t=0\},
\end{align}
where $T(0)=S(0)=0$. 
For $p \in [0 , 1]$, we define i.i.d. random variables ${\{ A^{[n]}_p \}}_{n\in \bN}$ as 
\begin{align}
A^{[n]}_p = 
\begin{cases}
0, ~~&\text{with probability }1-p, \\
1, ~~&\text{with probability }p. 
\end{cases}
\end{align}
We write 
\begin{align}
K^{p}_0 = &\kappa^{a^\ast , -}_0 \land \min \{ T[n] >0 : A^{[n]}_p = 0\},  \\
T^{p}_0 = &\tau^-_0 \land \min \{ S[n] >0 : A^{[n]}_p = 0\}. 
\end{align}
Then, we have 
\begin{align}
\beta \bE_{a^\ast}\sbra{ e^{-q K^{1}_0} }=\un{\nu}(a^\ast)<1\text{ and }\beta\bE_{a^\ast}\sbra{ e^{-q K^{0}_0} }=\lim_{a\uparrow a^\ast}\un{\nu}(a)\geq 1. 
\end{align}
By the strong Markov property, we have 
\begin{align}
\bE_{a^\ast}\sbra{ e^{-q K^{p}_0} }&= \bE_{a^\ast}\sbra{ e^{-q \kappa^{a^\ast , -}_0}; \kappa^{a^\ast , -}_0 <T[1] } 
+\bE_{a^\ast}\sbra{ e^{-q T[1]};  T[1]<\kappa^{a^\ast , -}_0 }(1-p)
\\
&+\bE_{a^\ast}\sbra{ e^{-q T[1]};  T[1]<\kappa^{a^\ast , -}_0 }
p\sum_{n=0}^\infty {\rbra{ \bE_{0}\sbra{ e^{-q T[1]};  T[1]<\kappa^{a^\ast , -}_0 } p} }^n \\
&\times\rbra{   \bE_0\sbra{ e^{-q T[1]};  T[1]<\kappa^{a^\ast , -}_0 }(1-p)         
+\bE_{0}\sbra{ e^{-q \kappa^{a^\ast {}, -}_0}; \kappa^{a^\ast , -}_0 <T[1] }        } , 
\end{align}
and so $p\mapsto \bE_{a^\ast}\sbra{ e^{-q K^{p}_0} }$ is continuous on $p \in [0, 1]$. 
Therefore, 
we can take some $p^\ast \in [0, 1)$ that satisfies \eqref{458}. 
\par
{\textbf{\textit{ii})}} 
In this step, we rewrite $v_{\pi^{a^\ast}}$ in a convenient form. 
For $x, a \in (0, \infty)$ with $x\leq a$ 
and $\epsilon>0$, we have 
\begin{align}
\un{\varphi}_{0, a} (x)
\leq 
\bE_x \sbra{e^{-q T^{p^\ast}_0}; T^{p^\ast}_0 < \tau^+_a }
\leq \un{\varphi}_{0,a }(x-\epsilon). 
\end{align}
By the continuity of $\un{\varphi}_{0,a }$, 
{for $a\in (0, \infty)$ and $x \in (0,a] \backslash E^{(1)}_a$, }
we have 
\begin{align}
\un{\varphi}_{0, a} (x)
= 
\bE_x \sbra{e^{-q T^{p^\ast}_0}; T^{p^\ast}_0 < \tau^+_a }. \label{460}
\end{align}
By the same argument as above, 
{for $a\in (0, \infty)$ and $x \in (0,a] \backslash E^{(1)}_a$, }
we have 
\begin{align}
\bar{\varphi}_{a, 0}(x)
=\bE_x\sbra{e^{-q\tau^+_a}; \tau^+_a < T^{p^\ast}_0}. \label{433}
\end{align}
So, by Lemma \ref{Thm203}, for 
{$a\in(0, \infty)$ and $x \in(0, a) \backslash E^{(1)}_a$, }
we have 
\begin{align}
v^{\prime}_{\pi^{a}}(x)=\bE_x\sbra{e^{-q\tau^+_a}; \tau^+_a < T^{p^\ast}_0}
+\beta \bE_x\sbra{e^{-qT^{p^\ast}_0};   T^{p^\ast}_0 < \tau^+_a}. \label{461}
\end{align} 
Because we have 
\begin{align}
\cbra{Y^{a^\ast}_t: t\in [0, \tau^+_{a^\ast})}=\cbra{X_t: t\in [0, \tau^+_{a^\ast})}, 
\end{align}
for $x\in (0,  a^\ast]$, we also have 
\begin{align}
\bE_x\sbra{e^{-qT^{p^\ast}_0};   T^{p^\ast}_0 < \tau^+_{a^\ast}}=\bE_x\sbra{e^{-qK^{p^\ast}_0} ; K^{p^\ast}_0< \tau^+_{a^\ast}},
\label{209}
\end{align}
and 
\begin{align}
\bE_{x}\sbra{e^{-qK^{p^\ast}_0};  \tau^+_{a^\ast}<K^{p^\ast}_0}
=&\bE_x \sbra{e^{-q\tau^+_{a^\ast}} ; \tau^+_{a^\ast}<K^{p^\ast}_0}
\bE_{a^\ast}\sbra{ e^{-qK^{p^\ast}_0}} \\
=&\bE_x \sbra{e^{-q\tau^+_{a^\ast}} ; \tau^+_{a^\ast}<T^{p^\ast}_0}
\bE_{a^\ast}\sbra{ e^{-qK^{p^\ast}_0}}. \label{211}
\end{align}
By \eqref{461}, \eqref{209}, and \eqref{211}, for $x\in (0, a^\ast)\backslash E^{(1)}_{a^\ast}$, we have 
\begin{align}
v_{\pi^{a^\ast}}^\prime (x)=
\frac{\bE_{x}\sbra{e^{-qK^{p^\ast}_0};  \tau^+_{a^\ast}<K^{p^\ast}_0}}{\bE_{a^\ast}\sbra{ e^{-qK^{p^\ast}_0}}}
+\beta \bE_x\sbra{e^{-qK^{p^\ast}_0} ; K^{p^\ast}_0< \tau^+_{a^\ast}}. \label{217}
\end{align}
By \eqref{458}, we have 
\begin{align}
\eqref{217}=\beta\bE_x \sbra{e^{-qK^{p^\ast}_0}}. \label{222}
\end{align}
From \eqref{458} and \eqref{222}, we obtain \eqref{419} and that $v_{\pi^{a^\ast}}$ is a concave function. 
\par
{Suppose that $X$ has unbounded variation paths. Then it is easy to check that 
the map $x \mapsto\beta \bE_{x}\sbra{e^{-qK^{p^\ast}_0}} $ is continuous. Thus we have}
\begin{align}
v_{\pi^{a^\ast}}^\prime (a^\ast-)&=
\beta \bE_{a^\ast}\sbra{e^{-qK^{p^\ast}_0}} =1  
\end{align}
and we obtain \eqref{420}. 
The proof is now complete. 
}
{
We define the Radon--Nikodym density $v_{\pi^{a^\ast}}^\prime$ by 
\begin{align}
v_{\pi^{a^\ast}}^\prime(x)= \beta\bE_x\sbra{e^{-q K^{p^\ast}_0}}, \quad x \in (0, a^\ast). \label{545}
\end{align}
In addition, when $X$ has unbounded variation paths, we define the Radon--Nikodym density $v_{\pi^{a^\ast}}^{\prime\prime}$ by
\begin{align}
v_{\pi^{a^\ast}}^{\prime\prime}(x)=
\bar{\varphi}_{a^\ast, 0}^\prime (x) +\beta \un{\varphi}_{0, a^\ast}^\prime(x), \quad &x \in (0, a^\ast). 
\end{align}
}
\begin{Lem}\label{Lem406}
The function $v_{\pi^{a^\ast}}$ belongs to $ C^{(1)}_{\text{line}}$. 
Furthermore, if $X$ has unbounded variation paths, then $v_{\pi^{a^\ast}} \in C^{(2)}_{\text{line}} $. 
\end{Lem}
\Proof{
By Assumption \ref{Ass101}, Lemma \ref{Thm203}, Lemma \ref{Lem404}, and Lemma \ref{Lem203}, 
it is obvious that $v_{\pi^{a^\ast}}$ belongs to $ C^{(1)}_{\text{line}}$. In addition, 
we know that 
$v_{\pi^{a^\ast}}$ {is continuously differentiable} 
{and $v_{\pi^{a^\ast}}^\prime$ has a Radon--Nikodym density on $(0, \infty)$ }when $X$ has unbounded variation paths. 
So, it is enough to check that $v_{\pi^{a^\ast}}^\prime$ has a locally bounded density on $(0, \infty)$ when $X$ has unbounded variation paths. 
\par 
Suppose $X$ has unbounded variation paths. 
We use the same notation as in the proof of Lemma \ref{Thm203}. 
For $x \in \bR$, let $Y^{(x)}$ be a reflected process defined by 
\begin{align}
Y^{(x)}_t=X^{(x)}_t -\rbra{ (\sup_{s \in [0, t]} X^{(x)}_s -a^\ast) \lor 0}, \quad t\geq 0, 
\end{align}
and write 
\begin{align}
\kappa^{(x)}= \inf \{t > 0 :  Y^{(x)}_t < 0 \}. 
\end{align}
We fix $a^\dagger \in (a^\ast , \infty)$. 
From {\eqref{545}} and because 
\begin{align}
 K^{p^\ast}_0 =K^{1}_0=\kappa^{a^\ast,-}_{0},
\end{align}
for $x \in (0, a^\ast) \backslash (E^{(2)}_{a^\ast}\cup E^{(2)}_{a^\dagger})$, we have 
\begin{align}
0\geq v_{\pi^{a^\ast}}^{\prime\prime}(x )= \lim_{\epsilon \downarrow 0}\frac{v_{\pi^{a^\ast}}^{\prime}(x+\epsilon ) -v_{\pi^{a^\ast}}^{\prime}(x )}{\epsilon}
=\lim_{\epsilon \downarrow 0}\frac{1}{\epsilon} \bE_0 \sbra{ e^{-q\kappa^{(x+\epsilon)}} -e^{-q\kappa^{(x)}}}. \label{546}
\end{align}
Because $Y^{(x+\epsilon)}_t = Y^{(x)}_t$ for $t \in [\tau^{(x),+}_{a^\ast}, \infty )$, we have 
\begin{align}
\kappa^{(x+\epsilon)}=\kappa^{(x)}, \quad \text{ on }\{ \tau^{(x),+}_{a^\ast} < \kappa^{(x)} \} ,
\end{align}
and so 
\begin{align}
\eqref{546}=&\lim_{\epsilon \downarrow 0}\frac{1}{\epsilon} \bE_0 \sbra{ e^{-q\kappa^{(x+\epsilon)}} -e^{-q\kappa^{(x)}} ; \kappa^{(x)}<\tau^{(x),+}_{a^\ast} }\\
=&\lim_{\epsilon \downarrow 0}\frac{1}{\epsilon} \bE_0 \sbra{ e^{-q\kappa^{(x+\epsilon)}} -e^{-q\kappa^{(x)}} ; \tau^{(x), -}_0<\tau^{(x),+}_{a^\ast} }. \label{549}
\end{align}
On $\{\tau^{(x), -}_0<\tau^{(x),+}_{a^\ast}\}$, we have 
\begin{align}
\kappa^{(x)}=\tau^{(x), -}_0 \leq \kappa^{(x+\epsilon)}\leq \tau^{(x+\epsilon), -}_0,
\end{align}
and so we have 
\begin{align}
\eqref{549} \geq&
\lim_{\epsilon \downarrow 0}\frac{1}{\epsilon} \bE_0 \sbra{ e^{-q\tau^{(x+\epsilon), -}_0} -e^{-q\tau^{(x), -}_0} ; \tau^{(x), -}_0<\tau^{(x),+}_{a^\ast} } \\ 
\geq&
\lim_{\epsilon \downarrow 0}\frac{1}{\epsilon} \rbra{\bE_0 \sbra{ e^{-q\tau^{(x+\epsilon), -}_0} ; \tau^{(x), -}_0<\tau^{(x),+}_{a^\ast} }
-\bE_0 \sbra{ e^{-q\tau^{(x), -}_0} ; \tau^{(x), -}_0<\tau^{(x),+}_{a^\ast} }}. \label{557}
\end{align}
Because 
\begin{align}
&\tau^{(x),-}_0\leq \tau^{(x+\epsilon),-}_0\text{ for all }\omega \in \Omega, \\
&\cbra{\tau^{(x), -}_0<\tau^{(x),+}_{a^\ast}}
\subset\cbra{\tau^{(x), -}_0<\tau^{(x),+}_{a^\dagger}}, \\
& \cbra{\tau^{(x+\epsilon), -}_0<\tau^{(x+\epsilon),+}_{a^\dagger}}
  \subset \cbra{\tau^{(x), -}_0<\tau^{(x),+}_{a^\dagger}},
\end{align}
we have
 \begin{align}
\eqref{557}\geq&
\lim_{\epsilon \downarrow 0}\frac{1}{\epsilon} \rbra{\bE_0 \sbra{ e^{-q\tau^{(x+\epsilon), -}_0} ; \tau^{(x+\epsilon), -}_0<\tau^{(x+\epsilon),+}_{a^\dagger} }
-\bE_0 \sbra{ e^{-q\tau^{(x), -}_0} ; \tau^{(x), -}_0<\tau^{(x),+}_{a^\dagger} }}\\
=&\lim_{\epsilon \downarrow 0}\frac{ \varphi_{0, a^\dagger} (x+\epsilon) -\varphi_{0, a^\dagger}(x)}{\epsilon}=
\varphi_{0, a^\dagger}^\prime (x). 
\end{align}
This implies that $v_{\pi^{a^\ast}}^{\prime\prime} $ is bounded on 
$(a^{\ddagger}, a^\ast){\backslash (E^{(2)}_{a^\ast}\cup E^{(2)}_{a^\dagger})}$ for $a^{\ddagger}\in (0, a^\ast)$. 
Because $v_{\pi^{a^\ast}}^{\prime\prime} = 0$ on ${[}a^\ast , \infty)$, 
the proof is complete. 
}
\begin{Lem}\label{Lem209}
Suppose $a^\ast > 0$. 
When $X$ has bounded variation paths, 
for $x \in (0, a^\ast)$, 
\begin{align}
\cL v_{\pi^{a^\ast}} (x) - q v_{\pi^{a^\ast}}(x) = 0.  \label{219}
\end{align}
{In addition, when $X$ has unbounded variation paths, we redefine $v_{\pi^{a^\ast}}^{\prime\prime}$ on $(0 , a^\ast)$ to satisfy \eqref{219}.}
\end{Lem}
\Proof{
For $x\in (0, a)$, we have 
\begin{align}
v^L_{\pi^{a^\ast}}(x) &= \bar{\varphi}_{a^\ast, 0}(x)
v^L_{\pi^{a^\ast}}(a^\ast)+\bE_x \sbra{  e^{-q\tau^+_{a^\ast}}(X_{\tau^+_{a^\ast}}-a^\ast) ;\tau^+_{a^\ast} <\tau^-_0  }
+\un{\varphi}_{0, a^\ast}(x)
v^L_{\pi^{a^\ast}}(0). 
\label{220}
\end{align}
For $x \in (0, a^\ast)$, the process $\cbra{M^{[1]}_t  :t\geq 0}$ where
\begin{align}
M^{[1]}_t =e^{-q(\tau^+_{a^\ast}\land \tau^-_0 \land t ) } 
\bar{\varphi}_{a^\ast, 0}\rbra{X_{\tau^+_{a^\ast} \land \tau^-_0 \land t }}
, \quad t\geq 0 ,
\end{align}
is a martingale under $\bP_x$ because 
\begin{align}
\bE_x \sbra{ e^{-q\tau^+_{a^\ast} }1_{\{ \tau^+_{a^\ast} < \tau^-_0\}} |\cF_t}
=&\bE_x \sbra{e^{-q\tau^+_{a^\ast} } 1_{\{\tau^+_{a^\ast} \leq t \land \tau^-_0\}}+e^{-q\tau^+_{a^\ast} } 1_{\{t< \tau^+_{a^\ast}<\tau^-_0 \}} |\cF_t}\\
=&e^{-q\tau^+_{a^\ast} } 1_{\{\tau^+_{a^\ast} \leq t \land \tau^-_0\}}+e^{-qt}  1_{\{ t<\tau^+_{a^\ast} \land \tau^-_0 \}}
\bar{\varphi}_{a^\ast, 0}\rbra{X_t}
\\
=&e^{-q(\tau^+_{a^\ast}\land \tau^-_0\land t ) } 
\bar{\varphi}_{a^\ast, 0}\rbra{X_{\tau^+_{a^\ast} \land \tau^-_0 \land t }}. 
\end{align}
By the same argument, for $x\in (0, a^\ast)$, 
$\cbra{M^{(2)}_t :t\geq 0}$ and $\cbra{M^{(3)}_t :t\geq 0}$ where
\begin{align}
M^{(2)}_t &=e^{-q(\tau^+_{a^\ast}\land\tau^-_{0}\land t ) } 
\un{\varphi}_{0, a^\ast}\rbra{X_{\tau^+_{a^\ast} \land \tau^-_0 \land t }}
,\quad &t\geq 0, \\
M^{(3)}_t &=e^{-q(\tau^+_{a^\ast}\land\tau^-_{0}\land t ) } \bE_{X_{\tau^+_{a^\ast} \land\tau^-_{0}\land t }} \sbra{e^{-q\tau^+_{a^\ast}} (X_{\tau^+_{a^\ast}}-a^\ast) ;\tau^+_{a^\ast} < \tau^-_0 }, \quad &t\geq 0, z
\end{align}
are martingales under $\bP_x$.  
By \eqref{220}, the process ${\cbra{ M^{(4)}_t:t\geq 0}} $ where 
\begin{align}
M^{(4)}_t =e^{-q(\tau^+_{a^\ast}\land\tau^-_{0}\land t ) } v_{\pi^a}^L(X_{\tau^+_{a^\ast}\land\tau^-_{0}\land t}), \quad t\geq 0, \label{568}
\end{align}
is a martingale under $\bP_x$. 
{By the argument used for the proof of \eqref{568}, the processes ${\cbra{ M^{(5)}_t:t\geq 0}} $ and ${\cbra{ M^{(6)}_t:t\geq 0}} $, where 
\begin{align}
M^{(5)}_t &=e^{-q(\tau^+_{a^\ast}\land\tau^-_{0}\land t ) } v_{\pi^a}^R(X_{\tau^+_{a^\ast}\land\tau^-_{0}\land t}), \quad t\geq 0, \\
M^{(6)}_t &=e^{-q(\tau^+_{a^\ast}\land\tau^-_{0}\land t ) } v_{\pi^a}(X_{\tau^+_{a^\ast}\land\tau^-_{0}\land t}), \quad t\geq 0, 
\end{align}
are martingales under $\bP_x$. }
By the same reasoning as that of the proof of \cite[(12)]{BifKyp2010}, for $x \in (0, a^\ast)\backslash E^{(1)}_{a^\ast}$ (resp., $ (0, a^\ast)\backslash 
E^{(2)}_{a^\ast}
$), we have \eqref{219}. 
{Here, we used the continuity of the map $x \mapsto \cL v_{\pi^{a^\ast}} (x)$ on 
$(0, a^\ast)\backslash E^{(1)}_{a^\ast}$ (resp., $ (0, a^\ast)\backslash 
E^{(2)}_{a^\ast}
$). }
\par
{When $X$ has bounded variation paths, we obtain \eqref{219} for $x \in (0, a^\ast)$ since 
Remark \ref{Rem204} holds and $v_{\pi^{a^\ast}}^\prime$ is right continuous. }
\par
{When $X$ has unbounded variation paths, by Remark \ref{Rem204}, we can redefine the locally bounded Radon--Nikodym density $v_{\pi^{a^\ast}}^{\prime\prime}$ on $(0, a^\ast)$, which is continuous almost everywhere and satisfies \eqref{219} for $x \in (0, a^\ast)$.}
\par
The proof is complete. 
}
\begin{Lem}\label{Lem408}
For $x \leq a^\ast$, we have 
\begin{align}
\cL v_{\pi^{a^\ast}} (x)-qv_{\pi^{a^\ast}}(x)\leq 0. \label{565}
\end{align}
\end{Lem}
\Proof{
This proof is almost the same as that of \cite[Lemma 5]{AvrPalPis2007}. 
We write $g(x)=\cL v_{\pi^{a^\ast}} (x)-qv_{\pi^{a^\ast}}(x)$ for $x\geq 0$. 
\par
From the form of the operator $\cL$, {\eqref{520}, and \eqref{544}}, for $x{\geq }a^\ast$, we have 
\begin{align}
g(x)=
 \gamma +\int_{\bR\backslash \{0\}}(v_{\pi^{a^\ast}}(x+z)-(x+b)-z1_{\{\absol{z}<1\}})\Pi(dz)- q(x+b), \label{294}
\end{align}
where $b= v_{\pi^{a^\ast}}(a^\ast) -a^\ast$. 
By the concavity of $v_{\pi^{a^\ast}}$ (see Lemma \ref{Lem203}) and the form of \eqref{294}, 
$g(x)$ is a continuous concave function on $(a^\ast , \infty)$. 
\par
We prove that, for $a>a^\ast$, 
\begin{align}
v_{\pi^a}(x)-v_{\pi^{a^\ast}}(x)
=\bE_x \sbra{\int_0^\infty e^{-qt}
g(U^{\pi^a}_t)1_{[a^\ast, \infty)}(U^{\pi^a}_t)dt}. \label{437}
\end{align}
We write $v_{\pi^{a^\ast}}^{(\epsilon)}(x)=v_{\pi^{a^\ast}}(x+\epsilon)$ for $x\in\bR$. 
Then, we can define $\cL v_{\pi^{a^\ast}}^{(\epsilon)}(x):=\cL v_{\pi^{a^\ast}}(x+\epsilon)$ for $x >-\epsilon$. 
Let $L^{\pi, c}$ be the continuous part of $L^\pi$ and let $R^{\pi, c}$ be the continuous part of $R^\pi$ for $\pi\in\cA$. 
By an application of the Meyer--It\^o formula (see \cite[Theorem IV.70 or IV.71]{Pro2005}), we have 
\begin{align}
\begin{aligned}
e^{-qt}v_{\pi^{a^\ast}}^{(\epsilon)}(U^{\pi^a}_t)-v_{\pi^{a^\ast}}^{(\epsilon)}(U^{\pi^a}_{0-})
=&-q\int_0^t e^{-qs}  v_{\pi^{a^\ast}}^{(\epsilon)}(U^{\pi^a}_{s-}) ds
+\int_0^t e^{-qs} v_{\pi^{a^\ast}}^{(\epsilon)\prime} (U^{\pi^a}_{s-})dU^{\pi^a}_s \\
&+\frac{\sigma^2}{2}\int_0^te^{-qs}v_{\pi^{a^\ast}}^{(\epsilon)\prime\prime}(U^{\pi^a}_{s-})ds \\
&+\sum_{0\leq s\leq t} e^{-qs}\rbra{ v_{\pi^{a^\ast}}^{(\epsilon)}(U^{\pi^a}_{s-}+\Delta U^{\pi^a}_s)-v_{\pi^{a^\ast}}^{(\epsilon)}(U^{\pi^a}_{s-}) -v_{\pi^{a^\ast}}^{(\epsilon)\prime} (U^{\pi^a}_{s-} ) \Delta U^{\pi^a}_s}. 
\end{aligned}
\label{569}
\end{align}  
Because 
\begin{align}
U^{\pi^a}_t = X_t -L^{\pi^a, c}_t - \sum_{0\leq s\leq t}\Delta L^{\pi^a}_s+R^{\pi^a, c}_t + \sum_{0\leq s\leq t}\Delta R^{\pi^a}_s, \quad t\geq 0, 
\end{align}
we have 
\begin{align}
\eqref{569}=&-q\int_0^t e^{-qs}  v_{\pi^{a^\ast}}^{(\epsilon)}(U^{\pi^a}_{s-}) ds
+\int_0^t e^{-qs} v_{\pi^{a^\ast}}^{(\epsilon)\prime} (U^{\pi^a}_{s-})dX_s 
-\int_0^t e^{-qs}v_{\pi^{a^\ast}}^{(\epsilon)\prime} (U^{\pi^a}_{s-})dL^{\pi,c}_s\\
&+\int_0^t e^{-qs}  v_{\pi^{a^\ast}}^{(\epsilon)\prime} (U^{\pi^a}_{s-})dR^{\pi,c}_s
+\frac{\sigma^2}{2}\int_0^te^{-qs}v_{\pi^{a^\ast}}^{(\epsilon)\prime\prime}(U^{\pi^a}_{s-})ds \\
&+\sum_{0\leq s\leq t} e^{-qs}\rbra{ v_{\pi^{a^\ast}}^{(\epsilon)}(U^{\pi^a}_{s-}+\Delta X_s)-v_{\pi^{a^\ast}}^{(\epsilon)}(U^{\pi^a}_{s-}) -v_{\pi^{a^\ast}}^{(\epsilon)\prime} (U^{\pi^a}_{s-} ) \Delta X_s}\\
&-\sum_{0\leq s\leq t} e^{-qs} \rbra{v_{\pi^{a^\ast}}^{(\epsilon)}(U^{\pi^a}_{s}+\Delta L^\pi_s)-v_{\pi^{a^\ast}}^{(\epsilon)}(U^{\pi^a}_{s})} \\
&-\sum_{0\leq s\leq t} e^{-qs} \rbra{v_{\pi^{a^\ast}}^{(\epsilon)}(U^{\pi^a}_{s}-\Delta R^\pi_s)-v_{\pi^{a^\ast}}^{(\epsilon)}(U^{\pi^a}_{s})}. 
\end{align}   
Rewriting the above equation leads to 
\begin{align}
\begin{aligned}
e^{-qt}v_{\pi^{a^\ast}}^{(\epsilon)}(U^{\pi^a}_t)-v_{\pi^{a^\ast}}^{(\epsilon)}(U^{\pi^a}_{0-})
&=\int_0^t e^{-qs}(\cL - q)v_{\pi^{a^\ast}}^{(\epsilon)}(U^{\pi^a}_{s-})ds
+M^{\pi^a}_t\\
&-\int_0^t e^{-qs}v_{\pi^{a^\ast}}^{(\epsilon)\prime} (U^{\pi^a}_{s-})dL^{\pi^a,c}_s
+\int_0^t e^{-qs}v_{\pi^{a^\ast}}^{(\epsilon)\prime} (U^{\pi^a}_{s-})dR^{\pi^a,c}_s\\
&-\sum_{0\leq s\leq t} e^{-qs} \rbra{v_{\pi^{a^\ast}}^{(\epsilon)}(U^{\pi^a}_{s}+\Delta L^{\pi^a}_s)-v_{\pi^{a^\ast}}^{(\epsilon)}(U^{\pi^a}_{s})}  \\
&-\sum_{0\leq s\leq t} e^{-qs} \rbra{v_{\pi^{a^\ast}}^{(\epsilon)}(U^{\pi^a}_{s}-\Delta R^{\pi^a}_s)-v_{\pi^{a^\ast}}^{(\epsilon)}(U^{\pi^a}_{s})}.
\end{aligned}
\label{575}
\end{align}
Here, $\{M^{\pi^a}_t: t\geq 0\}$ is a local martingale such that 
\begin{align}
M^{\pi^a}_t&=\sigma \int_0^t e^{-qs} v_{\pi^{a^\ast}}^{(\epsilon)\prime} (U^{\pi^a}_{s-})dB_s
\\
&+\int_{[0,t ] \times \bR}e^{-qs}
\rbra{ v_{\pi^{a^\ast}}^{(\epsilon)}(U^{\pi^a}_{s-}+y) - v_{\pi^{a^\ast}}^{(\epsilon)}(U^{\pi^a}_{s-}) 
 }(\cN (ds \times dy ) -ds \times \Pi(dy)) ,
\end{align}
where $B$ is a standard Brownian motion. 
Because $M^{\pi^a}$ is a local martingale, we can take a sequence of stopping times ${\{T^{\pi^a}_n\}}_{n\in\bN}$, which is a localizing sequence for $M$ with $T_n \uparrow \infty$ almost surely. 
We take the expectation of \eqref{575} at time $t\land T^{\pi^a}_n$ and 
take the limit as $t\uparrow \infty$ and $n\uparrow\infty$. By Lemma \ref{Lem408}, we have 
\begin{align}
-v_{\pi^{a^\ast}}^{(\epsilon)}(x)
&=\bE_x\sbra{ \int_0^\infty e^{-qs}(\cL - q)v_{\pi^{a^\ast}}^{(\epsilon)}(U^{\pi^a}_{s-}) 1_{[a^\ast-\epsilon , \infty)} (U^{\pi^a}_{s-})ds}\\
&-\bE_x\sbra{\int_0^\infty e^{-qs}v_{\pi^{a^\ast}}^{(\epsilon)\prime} (a)dL^{\pi^a,c}_s}
+\bE_x \sbra{\int_0^\infty e^{-qs}v_{\pi^{a^\ast}}^{(\epsilon)\prime} (0)dR^{\pi^a,c}_s}\\
&-\bE_x \sbra{\sum_{0\leq s\leq \infty} e^{-qs} \rbra{v_{\pi^{a^\ast}}^{(\epsilon)}(U^{\pi^a}_{s}+\Delta L^{\pi^a}_s)-v_{\pi^{a^\ast}}^{(\epsilon)}(U^{\pi^a}_{s})}  }\\
&-\bE_x\sbra{\sum_{0\leq s\leq \infty} e^{-qs} \rbra{v_{\pi^{a^\ast}}^{(\epsilon)}(U^{\pi^a}_{s}-\Delta R^{\pi^a}_s)-v_{\pi^{a^\ast}}^{(\epsilon)}(U^{\pi^a}_{s})}}.
\end{align}
By 
{\eqref{545}}
, $v_{\pi^{a^\ast}}^{\prime} (0+)=\beta$ when $0$ is regular for $(-\infty, 0)$. 
By contrast, 
$R^{\pi^a,c} \equiv 0$ when $0$ is irregular for $(-\infty, 0)$. So, we have 
\begin{align}
\bE_x \sbra{\int_0^\infty e^{-qs}v_{\pi^{a^\ast}}^{\prime} (0+)dR^{\pi^a,c}_s}
=\beta \bE_x \sbra{\int_0^\infty e^{-qs}dR^{\pi^a,c}_s}. \label{448}
\end{align} 
By taking the limit as $\epsilon \downarrow 0$ and applying Lemma \ref{Lem203}, \eqref{448}, and the continuity of $\cL v_{\pi^{a^\ast}}$, we have \eqref{437}. 
\par
By Lemma \ref{Lem202} and \eqref{437}, we have 
\begin{align}
\bE_x \sbra{\int_0^\infty e^{-qt}g(U^{\pi^a}_t)1_{[a^\ast, \infty)}(U^{\pi^a}_t)dt} \leq 0, \quad a   \in (a^\ast, \infty ).\label{449}
\end{align}
By the continuity and the concavity of $g$ along with \eqref{449}, we have \eqref{565} for $x \in [a^\ast , \infty)$. 
The proof is complete. 
}
\Proof[Proof of Theorem \ref{Thm501}]{
By Lemmas \ref{Lem203}, \ref{Lem406}, \ref{Lem209}, {\ref{Lem408}, and \eqref{545}}
$v_{\pi^{a^\ast}}\in C^{(1)}_{\text{line}}$ (resp., $C^{(2)}_{\text{line}}$) satisfies \eqref{201} and  \eqref{202} when $X$ has bounded (resp., unbounded) variation paths. 
From Proposition \ref{Prop201}, the proof is complete. 
}


\section{Examples}\label{Sec06}
{We assumed the continuity of $\bar{\varphi}_{a, 0}$ and $\un{\varphi}_{0, a}$ in Assumption \ref{Ass101} when $X$ has unbounded variation paths. }
In this section, we present examples of L\'evy processes {having unbounded variation paths} that satisfy Assumption \ref{Ass101}. 
\par
Let $X$ be a L\'evy process with characteristic exponent \eqref{202a} having unbounded variation paths. 
We additionally assume that $\Pi (-\infty, 0)< \infty$ or $\Pi ( 0, \infty)< \infty$. 
Then, for $a>0$, both $\bar{\varphi}_{a, 0}$ and $\un{\varphi}_{0, a}$ are continuously differentiable on $(0,a)$. 
We check this fact. 
\par
We assume without loss of generality that $\Pi ( 0, \infty)< \infty$. 
Then there exist a spectrally negative L\'evy process $Z$ with unbounded variation paths, 
a Poisson process $N^{(r)}$ with rate $r>0$, 
and i.i.d. positive random variables ${\{ J_n\}}_{n\in \bN}$ such that 
\begin{align}
X_t=Z_t + \sum_{i=1}^{N^{(r)}_t}J_i , \quad t \geq 0 . 
\end{align}
Here, $Z$ has the Laplace exponent $\psi_Z$, which satisfies
\begin{align}
e^{t\psi_Z (\lambda)}=\bE^{Z}_0\sbra{e^{\lambda Z_t}}, \quad \lambda \geq 0 , ~t \geq 0,  
\end{align}
where $\bP^Z_x$ is the law of $Z$ when it starts at $x\in\bR$. 
Then, $\psi_Z$ takes the form
\begin{align}
\psi_Z (\lambda) = \gamma\lambda +\frac{1}{2}\sigma^2 \lambda^2 
+\int_{(-\infty , 0)} (e^{\lambda x} -1 - \lambda x 1_{\{ x>-1\}}) \Pi_Z(dx) , \quad \lambda \in [0, \infty),  
\end{align}
where $\Pi_Z (\cdot) = \Pi (\cdot \cap (-\infty , 0))$. 
\par
We recall the definition of scale functions and some properties of these functions. 
For $p\geq 0$, let $W_Z^{(p)}$ be the scale function of $Z$, which is the function from $\bR$ to $[0, \infty)$ 
such that $W_Z^{(p)} = 0$ on $(-\infty , 0)$ and $W_Z^{(p)}$ on $[0, \infty)$ is continuous, satisfying
\begin{align}
\int_0^\infty e^{-\lambda x} W_Z^{(p)} (x) dx = \frac{1}{\psi_Z (\lambda) -p}, \quad \lambda >\Phi_Z(p),  
\end{align}
where $\Phi_Z (p)= \inf\{ s\geq 0 : \psi_Z(s)>p \}$. 
For proofs of uniqueness, existence, and the basic facts listed below, see, for example, \cite[Section 8]{Kyp2014}.
For $a>0$, $x\in [0, a]$, and a non-negative measurable function $f$, we have 
\begin{align}
\bE^Z_x \sbra{e^{-p\tau^+_a}; \tau^+_a< \tau^-_0}&=\frac{W_Z^{(p)} (x)}{W_Z^{(p)} (a)},\label{627} \\
\bE^Z_x \sbra{\int_0^{\tau^+_a \land \tau^-_0}\!e^{-pt} f(Z_t)dt }&=
\int_0^a \! f\rbra{y} \rbra{\frac{W_Z^{(p)}(x)}{W_Z^{(p)}(a)}W_Z^{(p)}(a-y)-W_Z^{(p)}(x-y)} dy  . 
\label{628}
\end{align}
The function $W_Z^{(p)}$ is continuously differentiable on $(0, \infty)$ 
because $Z$ has unbounded variation paths. 
\par
We prove that $\bar{\varphi}_{a, 0}$ is continuously differentiable. 
By the strong Markov property, we have, for $x \in (0, a)$,  
\begin{align}
\begin{aligned}
\bar{\varphi}_{a, 0} (x) = &\bE_x\sbra{e^{-q\tau^+_a}; \tau^+_a < \tau^-_0 \land T^{(r)}[1]}  \\
&+\bE_x \sbra{e^{-qT^{(r)}[1]}1_{\{ T^{(r)}[1]<\tau^+_a \land \tau^-_0\}} 
\bar{\varphi}_{a, 0} (X_{T^{(r)}[1]-} + J_1)}, 
\end{aligned}\label{630}
\end{align}
where $T^{(r)}[1]$ is the first jump time of $N^{(r)}$. By the definitions of $Z$ and $T^{(r)}[1]$, we have 
\begin{align}
\eqref{630}= \bE^Z_x \sbra{e^{-(q+r)\tau^+_a} ;\tau^+_a < \tau^-_0 }
+r \bE^Z_x \sbra{\int_0^{\tau^+_a \land \tau^-_0}e^{-(q+r)t} \bar{\varphi}_{a, 0}(Z_t+J_1)dt }. \label{631}
\end{align}
By \eqref{627} and \eqref{628}, we have 
\begin{align}
\eqref{631}= \frac{W_Z^{(q+r)} (x)}{W_Z^{(q+r)} (a)}+r 
\int_0^a \bE \sbra{ \bar{\varphi}_{a, 0}(y+J_1)} \rbra{\frac{W_Z^{(q+r)}(x)}{W_Z^{(q+r)}(a)}W_Z^{(q+r)}(a-y)-W_Z^{(q+r)}(x-y)} dy. 
\label{632}
\end{align}
We know that $\frac{W_Z^{(q+r)} (x)}{W_Z^{(q+r)} (a)}$ is continuously differentiable on $(0, a)$, so we consider 
the differentiability of 
\begin{align}
\int_0^a \bE \sbra{ \bar{\varphi}_{a, 0}(y+J_1)} \rbra{\frac{W_Z^{(q+r)}(x)}{W_Z^{(q+r)}(a)}W_Z^{(q+r)}(a-y)-W_Z^{(q+r)}(x-y)} dy. 
\label{635}
\end{align}
Here, we cannot obtain the derivative of \eqref{635} right away using the dominated convergence theorem 
because the derivative of $W_Z^{(q+r)\prime}$ may not be bounded. 
By \eqref{632} and Fubini's theorem, we have 
\begin{align}
\bar{\varphi}_{a, 0}(x)&=\frac{1}{W_Z^{(q+r)} (a)} \int_0^xW_Z^{(q+r)\prime} (z)dz\\
&+r 
\int_0^a \bE \sbra{ \bar{\varphi}_{a, 0}(y+J_1)} \rbra{\frac{W_Z^{(q+r)}(a-y)}{W_Z^{(q+r)}(a)} \int_0^xW_Z^{(q+r)\prime}(z)dz
-\int_0^xW_Z^{(q+r)\prime}(z-y)dz} dy
\\
&=\int_0^x h(z)dz, 
\end{align}
where 
\begin{align}
h(z)=
\begin{cases}
\frac{W_Z^{(q+r)\prime} (z)}{W_Z^{(q+r)} (a)}+r 
\int_0^a \bE \sbra{ \bar{\varphi}_{a, 0}(y+J_1)} 
\\
~~\times
\rbra{\frac{W_Z^{(q+r)}(a-y)}{W_Z^{(q+r)}(a)} W_Z^{(q+r)\prime}(z)
-W_Z^{(q+r)\prime}(z-y)} dy,\quad &z \in (0, a), \\
0, \quad &z \in (\infty , 0] \cup [a, \infty). 
\end{cases}
\end{align}
Here, $h \geq 0$ almost everywhere with respect to the Lebesgue measure because $\bar{\varphi}_{a, 0}$ is non-decreasing. 
So, by Fubini's theorem, we have 
\begin{align}
\eqref{635}&=
\int_0^a \bE \sbra{\int_0^{y+J_1} h(z)dz} \rbra{\frac{W_Z^{(q+r)}(x)}{W_Z^{(q+r)}(a)}W_Z^{(q+r)}(a-y)-W_Z^{(q+r)}(x-y)} dy
\label{640}\\
&=\bE \sbra{\int_{0}^{a+J_1} dz  h(z) \int_{0\lor {z-J_1}}^a\rbra{\frac{W_Z^{(q+r)}(x)}{W_Z^{(q+r)}(a)}W_Z^{(q+r)}(a-y)-W_Z^{(q+r)}(x-y)} dy}\\
&=\bE \sbra{\int_{0}^{a+J_1} dz  h(z) \rbra{\frac{W_Z^{(q+r)}(x)}{W_Z^{(q+r)}(a)}\int_{0\lor z-J_1}^a W_Z^{(q+r)}(a-y)dy -
\int_0^{x \land(x-z+J_1)}W_Z^{(q+r)}(y)dy}} . \label{642}
\end{align}
Because we have 
\begin{align}
\int_0^\infty h(z) dz \leq 1,
\end{align}
applying the dominated convergence theorem gives the derivative of \eqref{642} as 
\begin{align}
\bE \sbra{\int_{0}^{a+J_1} dz  h(z) \rbra{\frac{W_Z^{(q+r)\prime}(x)}{W_Z^{(q+r)}(a)}\int_{0\lor z-J_1}^a W_Z^{(q+r)}(a-y)dy -
W_Z^{(q+r)}(x \land(x-z+J_1))}},
\end{align}
which is continuous on $(0, a)$. 
Therefore, $\bar{\varphi}_{a, 0}$ is continuously differentiable on $(0, a)$.

\par
We can prove that $\un{\varphi}_{0, a}$ is continuously differentiable on $(0, a)$ by the same argument as above.

\appendix
\section{The behavior of $U^{\pi^a}$ and $U^{\pi^{a+\epsilon}}$}\label{SecA}
For $t \in [\un{\rho}^{a, [k-1]}_0, \bar{\rho}^{a, [k]}_a)$, the processes $U^{\pi^a}$ and $U^{\pi^{a+\epsilon}}$ behave as follows. 
We have 
\begin{align}
0 \leq U^{\pi^{a+\epsilon}}_{\un{\rho}^{a, [k-1]}_0-}-U^{\pi^a}_{\un{\rho}^{a, [k-1]}_0-}\leq \epsilon.  \label{411a}
\end{align} 
By the definitions of $\pi^a$ and $\bar{\rho}^{a, [k]}_a$, 
for $t \in [\un{\rho}^{a, [k-1]}_0 ,\bar{\rho}^{a, [k]}_a )$, processes $U^{\pi^a}$, $L^{\pi^a}$, and $R^{\pi^a}$ satisfy 
\begin{align}
&R^{\pi^a}_t =R^{\pi^a}_{\un{\rho}^{a, [k-1]}_0-}-\inf_{s\in[\un{\rho}^{a, [k-1]}_0, t]} ((U^{\pi^a}_{\un{\rho}^{a, [k-1]}_0-}+X_s -X_{\un{\rho}^{a, [k-1]}_0-})\land 0), 
\label{412b}\\
&U^{\pi^a}_t =U^{\pi^a}_{\un{\rho}^{a, [k-1]}_0-}+ (X_t -X_{\un{\rho}^{a, [k-1]}_0-})+R^{\pi^a}_t -R^{\pi^a}_{\un{\rho}^{a, [k-1]}_0-} \leq a , \label{412a}\\
&L^{\pi^a}_t =L^{\pi^a}_{\un{\rho}^{a, [k-1]}_0-}.  \label{414a}
\end{align}
Additionally, by the definition of $\pi^{a+\epsilon}$, processes $U^{\pi^{a+\epsilon}}$, $L^{\pi^{a+\epsilon}}$, and $R^{\pi^{a+\epsilon}}$ satisfy
\begin{align}
&R^{\pi^{a+\epsilon}}_t =R^{\pi^{a+\epsilon}}_{\un{\rho}^{a, [k-1]}_0-}- 
\inf_{s\in[\un{\rho}^{a, [k-1]}_0, t]} ((U^{\pi^{a+\epsilon}}_{\un{\rho}^{a, [k-1]}_0-}+X_s -X_{\un{\rho}^{a, [k-1]}_0-}) \land 0) , \label{414b}\\
&U^{\pi^{a+\epsilon}}_t =U^{\pi^{a+\epsilon}}_{\un{\rho}^{a, [k-1]}_0-}+ (X_t -X_{\un{\rho}^{a, [k-1]}_0-})+R^{\pi^{a+\epsilon}}_t -R^{\pi^{a+\epsilon}}_{\un{\rho}^{a, [k-1]}_0-} , \label{413a}\\
&L^{\pi^{a+\epsilon}}_t =L^{\pi^{a+\epsilon}}_{\un{\rho}^{a, [k-1]}_0-}  \label{A07}
\end{align}
before the right-hand side of \eqref{413a} hits $(a+\epsilon, \infty)$. 
From \eqref{411a}, \eqref{412b}, and \eqref{414b},
for $t \in [\un{\rho}^{a, [k-1]}_0 ,\bar{\rho}^{a, [k]}_a )$, 
\begin{align}
(\text{the right-hand side of } &\eqref{413a})-(\text{the right-hand side of } \eqref{412a}) \leq \epsilon. 
\end{align}
So the right-hand side of \eqref{413a} is no more than $a+\epsilon$ on $[\un{\rho}^{a, [k-1]}_0 ,\bar{\rho}^{a, [k]}_a )$, 
which implies that each of \eqref{414b}, \eqref{413a}, and \eqref{A07} holds for  $t \in [\un{\rho}^{a, [k-1]}_0 ,\bar{\rho}^{a, [k]}_a )$. 
From \eqref{411a}, \eqref{412b}, and \eqref{414b}, for $t \in [\un{\rho}^{a, [k-1]}_0 ,\bar{\rho}^{a, [k]}_a )$,
\begin{align}
\begin{aligned}
(R^{\pi^a}_t -R^{\pi^a}_{\un{\rho}^{a, [k-1]}_0-}) -&(R^{\pi^{a+\epsilon}}_t -R^{\pi^{a+\epsilon}}_{\un{\rho}^{a, [k-1]}_0-} )\\
&\in [0, U^{\pi^{a+\epsilon}}_{\un{\rho}^{a, [k-1]}_0-}-U^{\pi^a}_{\un{\rho}^{a, [k-1]}_0-}]   
 \text{ is non-decreasing} . 
\end{aligned}
\label{417a}
\end{align}
From \eqref{411a}, \eqref{412a}, \eqref{413a}, and \eqref{417a}, for $t \in [\un{\rho}^{a, [k-1]}_0 ,\bar{\rho}^{a, [k]}_a )$, 
\begin{align}
U^{\pi^{a+\epsilon}}_t -U^{\pi^{a}}_t  
\in [0, U^{\pi^{a+\epsilon}}_{\un{\rho}^{a, [k-1]}_0-}-U^{\pi^a}_{\un{\rho}^{a, [k-1]}_0-}]
\text{ is non-increasing}. \label{418a}
\end{align}
\par
For $t \in [\bar{\rho}^{a, [k]}_a, \un{\rho}^{a, [k]}_0)$, the processes $U^{\pi^a}$ and $U^{\pi^{a+\epsilon}}$ behave as follows. 
By the definitions of $\pi^a$ and $\un{\rho}^{a, [k]}_0$, for $t \in[\bar{\rho}^{a, [k]}_a ,\un{\rho}^{a, [k]}_0 )$, we have 
\begin{align}
&L^{\pi^a}_t =L^{\pi^a}_{\bar{\rho}^{a, [k]}_a-}+\sup_{s\in[\bar{\rho}^{a, [k]}_a ,t]} (U^{\pi^a}_{\bar{\rho}^{a, [k]}_a-}+X_s -X_{\bar{\rho}^{a, [k]}_a-}-a), \label{422a}\\
&U^{\pi^a}_t =U^{\pi^a}_{\bar{\rho}^{a, [k]}_a-}+ (X_t -X_{\bar{\rho}^{a, [k]}_a-})- (L^{\pi^a}_t -L^{\pi^a}_{\bar{\rho}^{a, [k]}_a-}) \geq 0, \label{421a}\\
&R^{\pi^a}_t =R^{\pi^a}_{\bar{\rho}^{a, [k]}_a-}. \label{426b}
\end{align}
Additionally, by the definition of $\pi^{a+\epsilon}$, the processes $U^{\pi^{a+\epsilon}}$, $L^{\pi^{a+\epsilon}}$, and $R^{\pi^{a+\epsilon}}$ satisfy
\begin{align}
&L^{\pi^{a+\epsilon}}_t =L^{\pi^{a+\epsilon}}_{\bar{\rho}^{a, [k]}_a-}+
\sup_{s\in[\bar{\rho}^{a, [k]}_a, t]} ((U^{\pi^{a+\epsilon}}_{\bar{\rho}^{a, [k]}_a-}+X_s -X_{\bar{\rho}^{a, [k]}_a-}-(a+\epsilon) )\lor 0),\label{424a}\\
&U^{\pi^{a+\epsilon}}_t =U^{\pi^{a+\epsilon}}_{\bar{\rho}^{a, [k]}_a-}+ (X_t -X_{\bar{\rho}^{a, [k]}_a-})- 
(L^{\pi^{a+\epsilon}}_t -L^{\pi^{a+\epsilon}}_{\bar{\rho}^{a, [k]}_a-}),  \label{423a}\\
&R^{\pi^{a+\epsilon}}_t =R^{\pi^{a+\epsilon}}_{\bar{\rho}^{a, [k]}_a-} \label{A16}
\end{align}
before the right-hand side of \eqref{423a} hits $(-\infty, 0)$. 
From \eqref{418a}, \eqref{422a}, and \eqref{424a}, 
for $t \in [\bar{\rho}^{a, [k]}_a ,\un{\rho}^{a, [k]}_0 )$, 
\begin{align}
(\text{the right-hand side of } &\eqref{423a})-(\text{the right-hand side of } \eqref{421a}) \geq 0. 
\end{align}
So the right-hand side of \eqref{423a} is non-negative on $[\bar{\rho}^{a, [k]}_a ,\un{\rho}^{a, [k]}_0 )$, 
which implies that each of \eqref{424a}, \eqref{423a}, and \eqref{A16} holds for $t \in[\bar{\rho}^{a, [k]}_a ,\un{\rho}^{a, [k]}_0 )$.
From \eqref{418a}, \eqref{422a}, and \eqref{424a}, for $t \in[\bar{\rho}^{a, [k]}_a ,\un{\rho}^{a, [k]}_0 )$, 
\begin{align}
\begin{aligned}
(L^{\pi^a}_t -L^{\pi^a}_{\bar{\rho}^{a, [k]}_a-}) -&(L^{\pi^{a+\epsilon}}_t -L^{\pi^{a+\epsilon}}_{\bar{\rho}^{a, [k]}_a-}) \\
&\in [0, \epsilon - (U^{\pi^{a+\epsilon}}_{\bar{\rho}^{a, [k]}_a-}-U^{\pi^a}_{\bar{\rho}^{a, [k]}_a-})]
\text{ is non-decreasing}. 
\end{aligned}\label{427a}
\end{align}
From \eqref{418a}, \eqref{421a}, \eqref{423a}, and \eqref{427a}, for $t \in[\bar{\rho}^{a, [k]}_a ,\un{\rho}^{a, [k]}_0 )$, 
\begin{align}
U^{\pi^{a+\epsilon}}_t  -U^{\pi^{a}}_t 
\in [  U^{\pi^{a+\epsilon}}_{\bar{\rho}^{a, [k]}_a-}-U^{\pi^a}_{\bar{\rho}^{a, [k]}_a-}, \epsilon]
\text{ is non-decreasing}. \label{425a}
\end{align}

\section{Proof of Proposition \ref{Prop201}} \label{Sec0B}
The proof of Proposition \ref{Prop201} is almost the same as that of \cite[Proposition 4 (ii)]{AvrPalPis2007}. 
\par
Let $\pi \in \Pi$ be any admissible strategy. Then, $U^\pi$ is a $[0, \infty)$-valued process. 
We fix $\epsilon >0 $ and define $w_\epsilon (x )= w (x + \epsilon)$ for $x\in\bR$. 
Then, we can define $\cL w_\epsilon (x )= \cL w (x + \epsilon)$ for $x>-\epsilon$. 
By \eqref{201} and \eqref{202}, we have
\begin{align}
\cL w_\epsilon (x)-q w_\epsilon (x) \leq 0,       ~~~~~~  &x \geq 0, \label{A01}\\ 
1 \leq w_\epsilon^\prime (x) \leq \beta, ~~~~~~&x \in \bR \backslash\{ c_{w_\epsilon}
\}. \label{A02}
\end{align}
By an application of the Meyer--It\^o formula (see \cite[Theorem II.31 and IV.71]{Pro2005})
and by the same calculation as used for \eqref{575}, 
we have 
\begin{align}
e^{-qt}w_\epsilon(U^\pi_t)-w_\epsilon(U^\pi_{0-})
&=\int_0^t e^{-qs}(\cL - q)w_\epsilon(U^\pi_{s-})ds
+M_t\\
&-\int_0^t e^{-qs}w_\epsilon^\prime (U^\pi_{s-})dL^{\pi,c}_s
+\int_0^t e^{-qs}w_\epsilon^\prime (U^\pi_{s-})dR^{\pi,c}_s  \label{A011}\\
&-\sum_{0\leq s\leq t} e^{-qs}\rbra{w_\epsilon(U^\pi_{s}+\Delta L^\pi_s)-w_\epsilon(U^\pi_{s})}  \label{A08}\\
&
-\sum_{0\leq s\leq t} e^{-qs}\rbra{w_\epsilon(U^\pi_{s}+\Delta L^\pi_s-\Delta R^\pi_s)-w_\epsilon(U^\pi_{s}+\Delta L^\pi_s)}
.\label{A09}
\end{align}
Here, $\{M_t : t\geq 0\}$ is a local martingale satisfying 
\begin{align}
M_t&=\sigma \int_0^t e^{-qs} w_\epsilon^\prime (U^\pi_{s-})dB_s
\\
&+\int_{[0,t ] \times \bR}e^{-qs}
\rbra{ w_\epsilon(U^\pi_{s-}+y) - w_\epsilon(U^\pi_{s-}) 
}(\cN (ds \times dy ) -ds \times\Pi(dy)) .
\end{align}
By \eqref{A02}, we have 
\begin{align}
&\eqref{A08}\leq -\int_0^t e^{-qs}dL^{\pi,c}_s
+\beta \int_0^t e^{-qs}dR^{\pi,c}_s, \\
&\eqref{A09}\leq -\sum_{0\leq s\leq t} e^{-qs} \Delta L^\pi_s, \quad \eqref{A09}\leq \beta \sum_{0\leq s\leq t} e^{-qs} \Delta R^\pi_s,
\end{align}
and so 
\begin{align}
e^{-qt}w_\epsilon(U^\pi_t)-w_\epsilon(U^\pi_{0-})
&\leq \int_0^t e^{-qs}(\cL - q)w_\epsilon(U^\pi_{s-})ds
+M_t-\int_{[0, t]} e^{-qs}dL^{\pi}_s
+\beta\int_{[0, t]} e^{-qs}dR^{\pi}_s. 
\end{align}
Because $M$ is a local martingale, we can take a sequence of stopping times ${\{T_n\}}_{n\in\bN}$ that is a localizing sequence for $M$ with $T_n \uparrow \infty$ almost surely. 
Then, taking an expectation, we have 
\begin{align}
w_\epsilon(x) &\geq\bE_x\sbra{ -\int_0^{t \land T_n} e^{-qs}(\cL - q)w_\epsilon(U^\pi_{s-})ds
+\int_{[0,{t \land T_n}]} e^{-qs}dL^{\pi}_s
-\beta\int_{[0, {t \land T_n}]} e^{-qs}dR^{\pi}_s} \\
&+\bE_x \sbra{e^{-q{(t\land T_n)}} w_\epsilon(U^{\pi}_{t\land T_n})}\\
&\geq\bE_x\sbra{ \int_{[0,{t \land T_n}]} e^{-qs}dL^{\pi}_s
-\beta\int_{[0,{t \land T_n}]} e^{-qs}dR^{\pi}_s} +w_\epsilon (0)\bE_x \sbra{e^{-q(t\land T_n)} }, \label{A15}
\end{align}
where in \eqref{A15} we used \eqref{A01}. 
By taking the the limit as $t\uparrow \infty$, $n \uparrow \infty$, and $\epsilon \downarrow 0$, the proof is complete.

\section{The behavior of $U^{(x)}$ and $U^{(x+\epsilon)}$ under $\bP_0$} \label{SecB}
In this section, we describe the behavior of $U^{(x)}$ and $U^{(x+\epsilon)}$ under $\bP_0$, which is necessary for the proof of Lemma \ref{Thm203}. 
We define the hitting times inductively as follows. For $n \in \bN$, 
\begin{align}
\un{\kappa}^{[0]}_0&=0, \\
\bar{\kappa}^{[n]}_a&=\inf \{ t\geq 
\un{\kappa}^{[n-1]} _0 : U^{(x+\epsilon)}_t =a\}, \\
\un{\kappa}^{[n]}_0&=\inf \{ t> \bar{\kappa}^{[n]}_a: U^{(x)}_t=0 \}. 
\end{align}
Then, we have the following by induction. \par
For $t \in [\un{\kappa}^{[n-1]}_0, \bar{\kappa}^{[n]}_a)$, the processes $U^{(x)}$ and $U^{(x+\epsilon)}$ behave as follows. 
We have 
\begin{align}
&0 \leq L^{(x+\epsilon)}_{\un{\kappa}^{[n-1]}_0-} -L^{(x)}_{\un{\kappa}^{[n-1]}_0-}\leq \epsilon ,\\
& 0 \leq R^{(x)}_{\un{\kappa}^{[n-1]}_0-}-R^{(x+\epsilon)}_{\un{\kappa}^{[n-1]}_0-}\leq \epsilon, \label{B08}\\
&
\begin{aligned}
& 0\leq  U^{(x+\epsilon)}_{\un{\kappa}^{[n-1]}_0-}-U^{(x)}_{\un{\kappa}^{[n-1]}_0-}  \\
&~=\epsilon-({ L^{(x+\epsilon)}_{\un{\kappa}^{[n-1]}_0-} -L^{(x)}_{\un{\kappa}^{[n-1]}_0-}})
-({ R^{(x)}_{\un{\kappa}^{[n-1]}_0-}-R^{(x+\epsilon)}_{\un{\kappa}^{[n-1]}_0-}}). 
\end{aligned}\label{B09}
\end{align}
\par
By the definitions of $\pi^a$ and $\bar{\kappa}^{[n]}_a$, for $t \in [\un{\kappa}^{[n-1]}_0, \bar{\kappa}^{[n]}_a)$, we have 
\begin{align}
&
L^{(x+\epsilon)}_t=L^{(x+\epsilon)}_{\un{\kappa}^{[n-1]}_0-}, \label{B10}\\
& R^{(x+\epsilon)}_t=R^{(x+\epsilon)}_{\un{\kappa}^{[n-1]}_0-}- \inf_{s \in  [\un{\kappa}^{[n-1]}_0, t]} 
(({  U^{(x+\epsilon)}_{\un{\kappa}^{[n-1]}_0-}+X_s -X_{\un{\kappa}^{[n-1]}_0-}  }) \land 0), \label{B11}
\\
&U^{(x+\epsilon)}_t=U^{(x+\epsilon)}_{\un{\kappa}^{[n-1]}_0-}  +(X_t -X_{\un{\kappa}^{[n-1]}_0-})+(R^{(x+\epsilon)}_t -R^{(x+\epsilon)}_{\un{\kappa}^{[n-1]}_0-}).  \label{B12}
\end{align}
Additionally, by the definition of $\pi^a$, the processes $U^{(x)}$, $L^{(x)}$, and $R^{(x)}$ satisfy
\begin{align} 
&
L^{(x)}_t=L^{(x)}_{\un{\kappa}^{[n-1]}_0-}, \label{B13}\\
& R^{(x)}_t=R^{(x)}_{\un{\kappa}^{[n-1]}_0-}- \inf_{s^\in  [\un{\kappa}^{[n-1]}_0, t]} 
({ ( U^{(x)}_{\un{\kappa}^{[n-1]}_0-}+X_s -X_{\un{\kappa}^{[n-1]}_0-} ) } \land 0), \label{B14}
\\
&U^{(x)}_t=U^{(x)}_{\un{\kappa}^{[n-1]}_0-} +(X_t -X_{\un{\kappa}^{[n-1]}_0-})+(R^{(x)}_t -R^{(x)}_{\un{\kappa}^{[n-1]}_0-}),   \label{B15}
\end{align} 
before the right-hand side of \eqref{B15} hits to $[a, \infty)$. 
By \eqref{B09}, \eqref{B11}, and \eqref{B14}, for $t \in [\un{\kappa}^{[n-1]}_0, \bar{\kappa}^{[n]}_a)$, 
\begin{align}
a>(\text{the right-hand side of \eqref{B12}})\geq (\text{the right-hand side of \eqref{B15}}), 
\end{align}
and so \eqref{B13}, \eqref{B14}, and \eqref{B15} hold for $t \in [\un{\kappa}^{[n-1]}_0, \bar{\kappa}^{[n]}_a)$. 
From \eqref{B08}, \eqref{B09}, \eqref{B11}, \eqref{B12}, \eqref{B14}, and \eqref{B15}, for $t \in [\un{\kappa}^{[n-1]}_0, \bar{\kappa}^{[n]}_a)$, we have 
\begin{align}
&R^{(x)}_t-R^{(x+\epsilon)}_t\in[0, \epsilon] \text{ is non-decreasing}; \label{B16}\\
&U^{(x+\epsilon)}_t-U^{(x)}_t\in[0, \epsilon] \text{ is non-increasing}. \label{B17}
\end{align}
From \eqref{B10}, \eqref{B13}, \eqref{B16}, and \eqref{B17}, we have 
\begin{align}
&0 \leq L^{(x+\epsilon)}_{\bar{\kappa}^{[n]}_a-} -L^{(x)}_{\bar{\kappa}^{[n]}_a-}\leq \epsilon ,\label{B18}\\
& 0 \leq R^{(x)}_{\bar{\kappa}^{[n]}_a-}-R^{(x+\epsilon)}_{\bar{\kappa}^{[n]}_a-}\leq \epsilon, \\
&
\begin{aligned}
& 0\leq  U^{(x+\epsilon)}_{\bar{\kappa}^{[n]}_a-}-U^{(x)}_{\bar{\kappa}^{[n]}_a-}  \\
&~=\epsilon-({ L^{(x+\epsilon)}_{\bar{\kappa}^{[n]}_a-} -L^{(x)}_{\bar{\kappa}^{[n]}_a-}})
-({ R^{(x)}_{\bar{\kappa}^{[n]}_a-}-R^{(x+\epsilon)}_{\bar{\kappa}^{[n]}_a-}}). 
\end{aligned}
\label{B20}
\end{align}
\par
For $t \in [\bar{\kappa}^{[n]}_a, \un{\kappa}^{[n]}_0)$, processes $U^{(x)}$ and $U^{(x+\epsilon)}$ behave as follows. 
By the definitions of $\pi^a$ and $\un{\kappa}^{[n]}_0$, for $t \in [\bar{\kappa}^{[n]}_a, \un{\kappa}^{[n]}_0)$, we have 
\begin{align}
&
L^{(x)}_t=L^{(x)}_{\bar{\kappa}^{[n]}_a-}+
\sup_{s\in [\bar{\kappa}^{[n]}_a, t]} 
(({  U^{(x)}_{\bar{\kappa}^{[n]}_a-}+X_s -X_{\bar{\kappa}^{[n]}_a-} -a })\lor 0), \label{B24}\\
& R^{(x)}_t=R^{(x)}_{\bar{\kappa}^{[n]}_a-}, \label{B25}
\\
&U^{(x)}_t=U^{(x)}_{\bar{\kappa}^{[n]}_a-} +(X_t- X_{\bar{\kappa}^{[n]}_a-})-(L^{(x)}_t- L^{(x)}_{\bar{\kappa}^{[n]}_a-}).  \label{B26}
\end{align} 
Additionally, by the definition of $\pi^a$, the processes $U^{(x+\epsilon)}$, $L^{(x+\epsilon)}$, and $R^{(x+\epsilon)}$ satisfy
\begin{align}
&
L^{(x+\epsilon)}_t=L^{(x+\epsilon)}_{\bar{\kappa}^{[n]}_a-}+
\sup_{s\in [\bar{\kappa}^{[n]}_a, t]} 
({  U^{(x+\epsilon)}_{\bar{\kappa}^{[n]}_a-}+X_s -X_{\bar{\kappa}^{[n]}_a-} -a }), \label{B21}\\
& R^{(x+\epsilon)}_t=R^{(x+\epsilon)}_{\bar{\kappa}^{[n]}_a-}, \label{B22}
\\
&U^{(x+\epsilon)}_t=U^{(x+\epsilon)}_{\bar{\kappa}^{[n]}_a-} +(X_t- X_{\bar{\kappa}^{[n]}_a-})-(L^{(x+\epsilon)}_t- L^{(x+\epsilon)}_{\bar{\kappa}^{[n]}_a-})  \label{B23}
\end{align}
before the right-hand side of \eqref{B23} hits $(-\infty, 0]$. 
From \eqref{B20}, \eqref{B24}, and \eqref{B21}, for $t \in [\bar{\kappa}^{[n]}_a, \un{\kappa}^{[n]}_0)$, 
\begin{align}
(\text{the right-hand side of \eqref{B23}})\geq (\text{the right-hand side of \eqref{B26}})>0, 
\end{align}
and so \eqref{B21}, \eqref{B22}, and \eqref{B23} hold for $t \in [\bar{\kappa}^{[n]}_a, \un{\kappa}^{[n]}_0)$. 
From \eqref{B18}, \eqref{B20}, \eqref{B24}, \eqref{B26}, \eqref{B21}, and \eqref{B23}, for $t \in [\bar{\kappa}^{[n]}_a, \un{\kappa}^{[n]}_0)$, we have 
\begin{align}
&L^{(x+\epsilon)}_t-L^{(x)}_t\in[0, \epsilon] \text{ is non-decreasing}; \label{B27}\\
&U^{(x+\epsilon)}_t-U^{(x)}_t\in[0, \epsilon] \text{ is non-increasing}. \label{B28}
\end{align}
From \eqref{B25}, \eqref{B22}, \eqref{B27}, and \eqref{B28}, we have 
\begin{align}
&0 \leq L^{(x+\epsilon)}_{\un{\kappa}^{[n]}_0-} -L^{(x)}_{\un{\kappa}^{[n]}_0-}\leq \epsilon ,\\
& 0 \leq R^{(x)}_{\un{\kappa}^{[n]}_0-}-R^{(x+\epsilon)}_{\un{\kappa}^{[n]}_0-}\leq \epsilon, \\
& 0\leq  U^{(x+\epsilon)}_{\un{\kappa}^{[n]}_0-}-U^{(x)}_{\un{\kappa}^{[n]}_0-} \n \\
&~=\epsilon-({ L^{(x+\epsilon)}_{\un{\kappa}^{[n]}_0-} -L^{(x)}_{\un{\kappa}^{[n]}_0-}})
-({ R^{(x)}_{\un{\kappa}^{[n]}_0-}-R^{(x+\epsilon)}_{\un{\kappa}^{[n]}_0-}})
. 
\end{align}


\section*{Acknowledgments}
I express my deepest gratitude to 
Prof. Kazutoshi Yamazaki and Prof. Kouji Yano 
for their comments and improvements. 
The author was supported by a KAKENHI grant, No. JP18J12680, from the Japan Society for the Promotion of Science.



\bibliographystyle{jplain}
\bibliography{NOBA_references_03} 

\begin{thebibliography}{10}

\bibitem{AvaSheBer2011}
B.~Avanzi, J.~Shen, and B.~Wong.
\newblock Optimal dividends and capital injections in the dual model with
  diffusion.
\newblock {\em Astin Bull.}, Vol.~41, No.~2, pp. 611--644, 2011.

\bibitem{AvrPalPis2007}
F.~Avram, Z.~Palmowski, and M.~R. Pistorius.
\newblock On the optimal dividend problem for a spectrally negative {L}\'evy
  process.
\newblock {\em Ann. Appl. Probab.}, Vol.~17, No.~1, pp. 156--180, 2007.

\bibitem{BauYam2015}
E.~J. Baurdoux and K.~Yamazaki.
\newblock Optimality of doubly reflected {L}\'{e}vy processes in singular
  control.
\newblock {\em Stochastic Process. Appl.}, Vol. 125, No.~7, pp. 2727--2751,
  2015.

\bibitem{BayKypYam2013}
E.~Bayraktar, A.~E. Kyprianou, and K.~Yamazaki.
\newblock On optimal dividends in the dual model.
\newblock {\em Astin Bull.}, Vol.~43, No.~3, pp. 359--373, 2013.

\bibitem{BifKyp2010}
E.~Biffis and A.~E. Kyprianou.
\newblock A note on scale functions and the time value of ruin for {L}\'evy
  insurance risk processes.
\newblock {\em Insurance Math. Econom.}, Vol.~46, No.~1, pp. 85--91, 2010.

\bibitem{BoSonTanWanYan2012}
L.~Bo, R.~Song, D.~Tang, Y.~Wang, and X.~Yang.
\newblock L\'{e}vy risk model with two-sided jumps and a barrier dividend
  strategy.
\newblock {\em Insurance Math. Econom.}, Vol.~50, No.~2, pp. 280--291, 2012.

\bibitem{CzaPerYam2018}
I.~Czarna, J.~L. P\'{e}rez, and K.~Yamazaki.
\newblock Optimality of multi-refraction control strategies in the dual model.
\newblock {\em Insurance Math. Econom.}, Vol.~83, pp. 148--160, 2018.

\bibitem{KuzKypRiv2012}
A.~Kuznetsov, A.~E. Kyprianou, and V.~Rivero.
\newblock The theory of scale functions for spectrally negative {L}\'evy
  processes.
\newblock In {\em L\'evy matters {II}}, Vol. 2061 of {\em Lecture Notes in
  Math.}, pp. 97--186. Springer, Heidelberg, 2012.

\bibitem{Kyp2014}
A.~E. Kyprianou.
\newblock {\em Fluctuations of {L}\'evy processes with applications}.
\newblock Universitext. Springer, Heidelberg, second edition, 2014.
\newblock Introductory lectures.

\bibitem{LiYin2019}
M.~Li and G.~Yin.
\newblock Optimal threshold strategies with capital injections in a spectrally
  negative {L}\'{e}vy risk model.
\newblock {\em J. Ind. Manag. Optim.}, Vol.~15, No.~2, pp. 517--535, 2019.

\bibitem{LiTanWanYan2014}
X.~Li, D.~Tang, Y.~Wang, and X.~Yang.
\newblock Optimal processing rate and buffer size of a jump-diffusion
  processing system.
\newblock {\em Ann. Oper. Res.}, Vol. 217, pp. 319--335, 2014.

\bibitem{NobPerYamYan2018}
K.~Noba, J.~L. P\'{e}rez, K.~Yamazaki, and K.~Yano.
\newblock On optimal periodic dividend and capital injection strategies for
  spectrally negative {L}\'{e}vy models.
\newblock {\em J. Appl. Probab.}, Vol.~55, No.~4, pp. 1272--1286, 2018.

\bibitem{PerYam2017_2}
J.~L. P\'erez and K.~Yamazaki.
\newblock On the optimality of periodic barrier strategies for a spectrally
  positive {L}\'evy process.
\newblock {\em Insurance Math. Econom.}, Vol.~77, pp. 1--13, 2017.

\bibitem{PerYamYu2018}
J.~L. P\'{e}rez, K.~Yamazaki, and X.~Yu.
\newblock On the bail-out optimal dividend problem.
\newblock {\em J. Optim. Theory Appl.}, Vol. 179, No.~2, pp. 553--568, 2018.

\bibitem{PerYam2017_1}
Jos\'{e}-Luis P\'{e}rez and Kazutoshi Yamazaki.
\newblock Refraction-reflection strategies in the dual model.
\newblock {\em Astin Bull.}, Vol.~47, No.~1, pp. 199--238, 2017.

\bibitem{Pro2005}
P.~E. Protter.
\newblock {\em Stochastic integration and differential equations}, Vol.~21 of
  {\em Stochastic Modelling and Applied Probability}.
\newblock Springer-Verlag, Berlin, 2005.
\newblock Second edition. Version 2.1, Corrected third printing.

\bibitem{Yam2018}
K.~Yamazaki.
\newblock Optimality of two-parameter strategies in stochastic control.
\newblock In {\em X{II} {S}ymposium of {P}robability and {S}tochastic
  {P}rocesses}, Vol.~73 of {\em Progr. Probab.}, pp. 51--104.
  Birkh\"{a}user/Springer, Cham, 2018.

\bibitem{YinSheWen2013}
C.~Yin, Y.~Shen, and Y.~Wen.
\newblock Exit problems for jump processes with applications to dividend
  problems.
\newblock {\em J. Comput. Appl. Math.}, Vol. 245, pp. 30--52, 2013.

\bibitem{YinYueShe2015}
C.~Yin, K.~C. Yuen, and Y.~Shen.
\newblock Convexity of ruin probability and optimal dividend strategies for a
  general {L}{\'e}vy process.
\newblock Vol. 2015, pp. 1--9, 2015.

\bibitem{YueYin2011}
K.~C. Yuen and C.~Yin.
\newblock On optimality of the barrier strategy for a general {L}\'{e}vy risk
  process.
\newblock {\em Math. Comput. Modelling}, Vol.~53, No. 9-10, pp. 1700--1707,
  2011.

\bibitem{ZhaCheYan2017}
Y.~Zhao, P.~Chen, and H.~Yang.
\newblock Optimal periodic dividend and capital injection problem for
  spectrally positive {L}\'{e}vy processes.
\newblock {\em Insurance Math. Econom.}, Vol.~74, pp. 135--146, 2017.

\bibitem{ZhaWanYaoChe2015}
Y.~Zhao, R.~Wang, D.~Yao, and P.~Chen.
\newblock Optimal dividends and capital injections in the dual model with a
  random time horizon.
\newblock {\em J. Optim. Theory Appl.}, Vol. 167, No.~1, pp. 272--295, 2015.

\bibitem{ZhaWanYin2017}
Y.~Zhao, R.~Wang, and C.~Yin.
\newblock Optimal dividends and capital injections for a spectrally positive
  {L}\'{e}vy process.
\newblock {\em J. Ind. Manag. Optim.}, Vol.~13, No.~1, pp. 1--21, 2017.

\end{thebibliography}

\end{document}